\documentclass[times, doublespace]{fldauth}
\usepackage{graphicx}
\usepackage{amssymb,amsxtra,amsmath,amsbsy}
\usepackage{latexsym}
\usepackage{epstopdf}

\usepackage{verbatim}
\usepackage{enumerate}
\usepackage{verbatim}
\usepackage{ifthen}
\usepackage{calc,fancyhdr}
\usepackage{indentfirst}
\usepackage{parskip}
\usepackage{cite}

\usepackage{color}

\usepackage{def}

\begin{document}


\title{{Radial Basis Function (RBF)-based} Parametric Models for Closed and Open Curves within the Method of Regularized Stokeslets}

\author{Varun Shankar\affil{1}\corrauth \ and Sarah D. Olson\affil{2}}

\address{\affilnum{1}Department of Mathematics, University of Utah, Salt Lake City, UT, USA\break
\affilnum{2}Department of Mathematical Sciences, Worcester Polytechnic Institute, Worcester, MA, USA}

\corraddr{Department of Mathematics, University of Utah, Salt Lake City, UT, USA. E-mail: vshankar@math.utah.edu}


\bigskip
\begin{abstract}
The method of regularized Stokeslets (MRS) is a numerical approach using regularized fundamental solutions to compute the flow due to an object in a viscous fluid where inertial effects can be neglected. The elastic object is represented as a Lagrangian structure, exerting point forces on the fluid. The forces on the structure are often determined by a bending or tension model, previously calculated using finite difference approximations. In this paper, we study Spherical Basis Function (SBF), Radial Basis Function (RBF) and Lagrange-Chebyshev parametric models to represent and calculate forces on elastic structures that can be represented by an open curve, motivated by the study of cilia and flagella. {The evaluation error for static open curves for the different interpolants, as well as errors for calculating normals and second derivatives using different types of clustered parametric nodes, are given for the case of an open planar curve.} We determine that SBF and RBF interpolants built on clustered nodes are competitive with Lagrange-Chebyshev interpolants for modeling twice-differentiable open planar curves. 
We propose using SBF and RBF parametric models within the MRS 
for evaluating and updating the elastic structure. Results for open and closed elastic structures immersed in a 2D fluid are presented, showing the efficacy of the {RBF}-Stokeslets method.
\end{abstract}

\maketitle


\section{Introduction}
\noindent
The method of regularized Stokeslets (MRS), developed by Cortez \cite{Cortez2000,Cortez05}, is a numerical method used to determine fluid flow in the presence of an elastic structure at zero Reynolds number. Using the linearity of the Stokes Equations, the velocity field is calculated as a superposition of regularized fundamental solutions in which the singularity has been removed.  Since microorganisms live in a viscosity dominated environment where the effects of inertia can be neglected, the MRS has been successfully used to model bacteria (e.g. \textit{E. coli} and spirochetes) \cite{Cisneros08,Cortez05,Flores05,Jung07}, dinoflagellates \cite{Nguyen11b}, sperm \cite{Gillies09,Olson11,Olson13a,Olson13b,Olson14a,Simons14}, cilia \cite{Ainley08,Leiderman13,Leiderman14}, biofilms \cite{Cogan05}, filaments \cite{Bouzarth11,Lee14}, and biflagellated algae \cite{OMalley12}. 

In each of these applications, the force that the organism exerts on the surrounding fluid is modeled by a regularized force term in the Stokes equation. 
We consider elastic structures treated as immersed boundaries, where {the curve corresponds to a centerline approximation of a slender body} and can be open or closed. The curves can be described in a Lagrangian form that can be discretized into a set of points whose trajectories will be calculated and updated using the MRS. These elastic structures generate bending or stretching forces via virtual springs or connections between the discretized points, which have been previously calculated by finite difference approximations. Alternative representations of structures via an interpolant and direct calculation of forces from interpolants has not been explored in the context of the MRS; however, this has been completed for closed curves in the Immersed Boundary (IB) method \cite{Shankar2012,SWKFJCP2014}. In this work, we are interested in exploring two variants of Radial Basis Function (RBF) interpolants to describe Lagrangian elastic structures within the MRS.

RBF interpolation is a powerful tool for the approximation of data at scattered node locations in multiple dimensions. On planar domains, RBF interpolants can be viewed as a generalization of polynomial interpolants (which are generally guaranteed to be unisolvent only on scattered nodes in 1D)~\cite{DriscollFornberg2002}. In this work, however, we are interested in the comparison of RBF interpolants restricted to a circle~\cite[\S 6]{FasshauerSchumaker:1998},~\cite[Ch. 17]{Wendland:2004} to more standard RBF interpolants. {Studies have shown that these restricted RBF interpolants (also known as Spherical Basis Functions or SBFs) provide spectral accuracy for smooth target functions~\cite{JetterStocklerWard:1999}, and error estimates for target functions in Sobolev spaces are also known~\cite{NarcSunWard:2007}. To date, a complete comparison of SBF interpolants and more standard RBF interpolants in the parametric interpolation setting has not been completed}. Of greater interest, SBF interpolants have been shown to be generalizations of Fourier-based interpolants~\cite{FornbergPiret:2007} and exhibit superconvergence in the approximation of derivatives of functions on circles, spheres and tori~\cite{FuselierWright2014}.

In previous work~\cite{Shankar2012}, we presented a parametric SBF representation of the boundaries of closed objects simulated in the IB method and showed that it was more accurate and less costly than the collection of techniques used in the traditional IB method (using piecewise quadratics to compute normal vectors and using finite-differences to compute forces). In addition, we showed that the parametric SBF model offered accuracy and convergence comparable to Fourier-based methods when the target shapes were infinitely smooth, and better accuracy and convergence than Fourier-based methods when the target shapes had only one or two underlying derivatives. In subsequent work~\cite{SWKFJCP2014}, we used the SBF geometric model within a full IB simulation of closed surfaces (platelets). The SBF geometric model has since been used within the Augmented Forcing method, a second-order forcing method for fluid-structure interaction~\cite{SWFKIJNMF2014}, and for forming surfaces out of point cloud data (on which PDEs were then solved)~\cite{FuselierWright2012, SWFKJSC2014}. However, in all these scenarios, the SBF geometric model was used only for the modeling of closed curves and/or surfaces, and the basis functions themselves were always restricted to the circle/sphere.

In this work, we explore the use of this SBF geometric model for the simulation of elastic structures, represented by both open and closed curves, within the MRS. Motivated by the application of modeling cilia and flagella, we first focus on a detailed comparison of SBF, RBF and Lagrange-Chebyshev geometric models for a Lagrangian structure corresponding to a planar, open curve. A static test case of a filament corresponding to a perturbed sinusoidal shape, a simplified centerline representation of a sperm flagellum,  is used to investigate interpolation error as well as error in computing normals and second derivatives. These studies are performed on Chebyshev nodes and the so-called ``mapped Chebyshev'' nodes in parametric space, to help determine which node sets are well-suited for the modeling of {perturbations of idealized shapes}. We then present results for time dependent fluid-structure interaction using an SBF parametric model within the MRS. We test the accuracy of the SBF and RBF models within this method for a test case with an exact solution in comparison to finite difference approximations of forces. Results of time-dependent fluid-structure interaction simulations involving closed planar shapes and an open filament propagating a sinusoidal wave are presented. The proposed RBF-Stokeslets method compares well with the MRS using finite difference approximations for the forces and {computation time could be potentially decreased with proper choice of nodes}.

The remainder of this article is organized as follows. In Section 2, we describe the adaptation of the SBF geometric model to the modeling of open and closed planar curves, and present a simple new RBF model for the same. We review the method of Regularized Stokeslets in Section 3, and describe how we use the {SBF and RBF models} within that method; specifically, we discuss modifications to the standard Stokeslet algorithm to incorporate the two RBF models. In Section 4, we present results for error on static test cases as well as time dependent simulation results for a dynamic elastic structure using the SBF-based regularized Stokeslets method. We conclude with a summary of our results, and outline possible future work in Section 5.

\section{Geometric Modeling with Radial Basis Functions}
\label{sec:geommodel}

In this section, we will briefly review the previously developed SBF model for closed planar curves and adapt the SBF representation to the modeling of open curves in 2D domains. We choose to use a parametric representation for these planar curves, since our target shapes (flagella, cilia, etc.) can naturally be represented in this fashion \cite{Gillies09,Olson11,Leiderman13}. 

We represent a general planar curve at any time $t$ parametrically by
\begin{equation}
\label{eq:2D_obj1}
\vX(\lon,t) = (X(\lon,t),Y(\lon,t)),
\end{equation}
where $\lon$ is the parametric variable. If these curves are open and planar (and of the form $Y = f(X)$ where $f$ is some smooth function), for a specific time $t$, we can simplify this parametrization to
\begin{equation}
\label{eq:2D_obj}
\vX(\lon,t) = (\lon, Y(\lon,t)).
\end{equation}
We explicitly track a finite set of $N_d$ points $\vXd_1(t),\ldots,\vXd_{N_d}(t)$, which we refer to as \emph{data sites}.  Here $\vXd_j(t) := \vX(\lond_j,t)$, $j=1,\ldots,N_d$, and we refer to the parameter values $\lond_1,\ldots,\lond_{N_d}$ as the \emph{data site nodes} (or simply \emph{nodes}). For a general planar curve, we construct each component of $\vX$ by using a smooth parametric RBF interpolant of each coordinate of the data sites.

{We now explain how to construct two types of RBF interpolants to the y-coordinate $Y(\lambda,t)$ using the data $(\lond_1,\Yd_1(t)),...,(\lond_{N_d},\Yd_{N_d}(t))$. Define $Y(\lon,t)$ by
\begin{align}
Y(\lon,t) &= \sum_{k=1}^{N_d} c^Y_k(t) \phi\lf(r_{\lambda,k}\rt),
\label{eq:circ_rbf_interp}
\end{align}
for coefficients $c^Y_k(t)$ and scalar-valued radial kernel $\phi(r)$, which we discuss below. 
The argument $r$ corresponds to the Euclidean distance between the points on the unit circle for a Spherical Basis Function (SBF) and a standard distance argument for an RBF interpolant. Specifically, we define $r$ as follows for the SBF and RBF interpolant,
\begin{subequations}
\begin{align}
r_{\lambda,k}&=\sqrt{2 - 2\cos(\lon-\lond_k)}\hspace{0.2cm}(SBF), \label{SBFr}\\
r_{\lambda,k}&=|\lon-\lond_k|\hspace{0.2cm}(RBF), \label{RBFr}
\end{align}
\end{subequations}
corresponding to the distance between a point $\lon$ and $\lon_k$.
For the geometric modeling of closed and open planar curves, we use the multiquadric (MQ) radial kernel function, given explicitly by
\begin{align}
\phi(r) = \sqrt{1 + (\ep r)^2} \label{eq:mq}, 
\end{align}
where $\ep > 0$ is called the shape parameter. Regardless of whether we use the SBF or the standard RBF, to
have $Y(\lon,t)$ interpolate the given data, we require that the
coefficients $c_k^Y(t), k=1,...,N_d$ for some $t$ satisfy the following system of
equations:}
\begin{align}
\underbrace{
\begin{bmatrix}
\phi\lf(r_{1,1}\rt) & \cdots & \phi\lf(r_{1,N_d}\rt) \\
\phi\lf(r_{2,1}\rt) & \cdots & \phi\lf(r_{2,N_d}\rt) \\
\vdots & \ddots & \vdots \\
\phi\lf(r_{N_d,1}\rt) & \cdots & \phi\lf(r_{N_d,N_d}\rt)
\end{bmatrix}}_{\displaystyle A}
\underbrace{
\begin{bmatrix}
c^Y_1(t) \\
c^Y_2(t) \\
\vdots \\
c^Y_{N_d}(t)
\end{bmatrix}}_{\displaystyle \ucdb}
= 
\underbrace{
\begin{bmatrix}
\Yd_1(t) \\
\Yd_2(t) \\
\vdots \\
\Yd_{N_d}(t)
\end{bmatrix}}_{\displaystyle \uYd },
\label{eq:rbf_linsys}
\end{align}
{where $r_{j,k}$ will be the distance between $\lond_j$ and $\lond_k$ using either Equation \eqref{SBFr} or \eqref{RBFr} for the SBF or RBF, respectively}.  Since $r_{j,k} = r_{k,j}$, the matrix $A$ is symmetric. For the MQ kernel in conjunction with the SBF, $A$ is non-singular if $\lond$ is not a multiple of $2\pi$; the standard RBF is non-singular for all distinct values of $\lond$~\cite{Fasshauer:2007,Wendland:2004}. Though $A$ may be highly ill-conditioned for small values of $\ep$, the interpolant is well-defined even in the limit $\ep \to 0$~\cite{DriscollFornberg2002,FornbergPiret:2007}. We will compare the performance of both SBFs and RBFs on a static test problem and comment on selecting $\ep$ in Section 4. The coefficient vector $\ucd$ can be computed using the same interpolation matrix $A$ and setting the right hand side of Equation \eqref{eq:rbf_linsys} to be the vector $\uXd$. We note that our choice of the multiquadric is somewhat arbitrary. One could use any of the infinitely-smooth RBFs in these scenarios. 

A somewhat more standard choice for this form of 1D interpolation problem would be to use a polynomial interpolant at Chebyshev nodes. Indeed, these interpolants are known to converge at an exponential/geometric rate when approximating analytic (or sometimes even simply smooth) functions. For completion, we therefore also present a well-known stable algorithm for computing Lagrange interpolants at Chebyshev nodes. In this work, we will refer to these interpolants as Lagrange-Chebyshev interpolants, since the term ``Chebyshev interpolant'' could also refer to interpolation with a Chebyshev polynomial (which we do not consider in this work).

Assuming a parametrization of open planar curves on an interval {$[a,b]$}, the corresponding Chebyshev nodes are given by
\begin{align}
\lond_k = \frac{b-a}{2} + \frac{b-a}{2}\cos\lf(\frac{2k-1}{2 N_d}\pi\rt),
\label{eq:cheb_nodes}
\end{align}
for $k=1,\ldots,N_d$. Note that the Chebyshev nodes do not include the endpoints of the interval. We must now form the polynomial interpolant at the Chebyshev nodes. It is well-known that forming the Vandermonde interpolation matrix corresponding to the monomials is extremely ill-conditioned at equispaced nodes. To ameliorate this difficulty, we form the Lagrange interpolating polynomial at the Chebyshev nodes. For the $Y(\lambda,t)$ component, this may be written as
{\begin{align}
Y(\lon,t) = \sum_{k=1}^{N_d}\Yd_k\ell_k(\lon),
\end{align}
where
\begin{align}
\ell_k(\lon) = \prod_{\substack{1 \leq m \leq k \\m \neq k}} \frac{\lon - \lond_m}{\lond_k - \lond_m},
\end{align}}are the Lagrange basis polynomials. To further ameliorate the ill-conditioning, we choose to use the \emph{barycentric form} of the Lagrange interpolating polynomial. First, we define the barycentric weights to be
{\begin{align}
w_k= \frac{1}{\prod_{\substack{1 \leq j \leq N_d \\ k\neq j} } (\lond_k - \lond_j)},
\end{align}}which allows us to write the Lagrange basis functions as $\ell_k(\lon) = \ell(\lon)\frac{w_k}{\lon - \lond_k}$, where $\ell(\lon) = (\lon - \lond_1)(\lon - \lond_2)\hdots(\lon-\lond_{N_d})$. The barycentric form of the Lagrange interpolant is very efficient to use for a fixed parametrization, and quite stable for hundreds of points.

{As mentioned earlier, we model \emph{open} planar curves by the simpler parametrization $\vX(\lon,t) = (\lon, Y(\lon,t))$. In this case, since the parametrization is a simple linear function in the x-direction, we only need to interpolate the y-coordinate of the data sites and compute the coefficient vector $\ucdb$. All derivatives of the x-coordinate can be computed analytically. If one does require an interpolant $X(\lon,t)$, we recommend replacing the scalar-valued radial kernel $\phi(r)$ in Equation \eqref{eq:circ_rbf_interp} with the term $|\lon - \lond_k|$. This is equivalent to a linear spline (which would be exact for linear functions). In our experiments, we also found that simply using the RBF in Equation \eqref{RBFr} also produces reasonable results, with the RBF interpolant being generally better conditioned than the linear spline. }

Once we have an {SBF or RBF interpolant}, we may wish to evaluate the interpolant at some set of $N_s$ \emph{evaluation nodes}, $\{\lons_j\}_{j=1}^{N_s}$, where $N_s >> N_d$. Further, we may wish to compute derivatives of the interpolant and evaluate those derivatives at either the data site nodes or the evaluation nodes. We will now describe how to perform these operations without explicitly computing the coefficient vectors $\ucd$ and $\ucdb$. Again, for simplicity, we focus on evaluating and differentiating the interpolant $Y(\lon,t)$. Formally, from Equation \eqref{eq:rbf_linsys}, we have
\begin{align}
\ucdb = A^{-1}\uYd.
\end{align}
We call the Cartesian points corresponding to the evaluation nodes the \emph{sample sites}. Now, let $\uYs$ be the vector with $N_s$ entries that contains the y-coordinates of the sample sites, obtained by evaluating the interpolant $Y(\lon,t)$ at the set of evaluation nodes. To obtain the vector $\uYs$, we perform the following operation:
\begin{align}
\uYs = 
B\ucdb = BA^{-1}\uYd \label{eq:nodeeval}, 
\end{align}
where $\{B_{ij}\}=\phi(r_{j,k})$ for $j=1,2,\hdots,N_s, k=1,2,\hdots,N_d$, and $r_{j,k}$ is defined by Equation \eqref{SBFr} or \eqref{RBFr} for the SBF and RBF, respectively.  ${E} = BA^{-1}$ is the time-independent $N_s \times N_d$ \emph{evaluation matrix}. One can obtain the x-coordinates of the sample sites similarly by the operation $\uXs = {E}\uXd$. Defining the matrices $\uvXd = [\uXd \ \uYd]$ and $\uvXs = [\uXs \ \uYs]$, then $\uvXs = {E}\uvXd$. For open planar curves, we find that it is sufficiently accurate to use a single evaluation matrix ${E}$, though one could also take the approach of defining two evaluation matrices (one corresponding to a piecewise-linear interpolant for the x-coordinate and the other corresponding to an SBF or RBF interpolant for the y-coordinate).

In order to compute quantities like tangents, normals and forces, we also need to be able to compute derivatives of the coordinate interpolants with respect to the parameter $\lon$ and evaluate them at either the data site nodes $\lond$ or the evaluation nodes $\lons$. To compute the $n^{th}$-order derivative of either an SBF or RBF interpolant and evaluate the derivative at the \emph{data site nodes}, we define the matrix
\begin{align}
\{B_{\lond}^n\}_{j,k} = \lf.\frac{\partial^n}{\partial \lon^n} \phi\lf(r(\lon) \rt) \rt|_{\lon = \lond_j}, \ j,k=1,\hdots,N_d,
\end{align}
which is an $N_d \times N_d$ matrix. To compute the $n^{th}$-order derivative of an RBF interpolant and evaluate the derivative at the \emph{evaluation nodes}, we similarly define the matrix $\{B_{\lons}^n\}_{j,k}$ where $j=1,\hdots,N_s$, and $k=1,\hdots,N_d$.
If we explicitly computed the coefficient vectors $\ucd$ and $\ucdb$, the products of the matrices $B_{\lond}^n$ and $B_{\lons}^n$ with those vectors would yield derivatives of interpolants at the data sites nodes and evaluation nodes respectively. However, since it is cheaper to use a coefficient-free formulation, we will now define analogs to Equation \eqref{eq:nodeeval} that perform coefficient-free differentiation. We thus define the following \emph{differentiation matrices}:
\begin{equation}
\mathcal{D}_{\lond}^n = B_{\lond}^n A^{-1}, \hspace{0.3cm}\mathcal{D}_{\lons}^n = B_{\lons}^n A^{-1},\label{diffM}
\end{equation}
which yield derivatives of the $X(\lon,t)$ and $Y(\lon,t)$ interpolants at data site nodes and sample site nodes when post-multiplied by the vectors $\uXd$ and $\uYd$. These differentiation matrices thus combine the operations of differentiation and evaluation of the interpolants. 

It is important to note that since the evaluation and differentiation matrices implicitly contain the interpolation operation (because they contain the term $A^{-1}$), these matrices can be applied to any $N_d$-long vector of function samples at the data sites to obtain either samples of the function at sample sites, or samples of derivatives of the function at data sites or sample sites. For a more detailed description of these interpolation and differentation matrices, we refer the reader to previous work~\cite{SWKFJCP2014}.

The barycentric form allows for fast evaluation of the Lagrange-Chebyshev interpolant. It also allows one to compute derivatives quite easily at the data site nodes. However, it is not clear how to evaluate the derivatives at the sample site nodes using this form of the interpolant, unlike with RBF/SBF or other interpolants. However, we use a simple alternative to evaluate derivatives.  We first use the barycentric form of the interpolant to obtain derivatives of the interpolant at the data site nodes; then, we once again use the barycentric form of the interpolant to interpolate the derivative samples themselves, and evaluate \emph{that} interpolant at the sample site nodes. For polynomials, we believe this process is quite reasonable, with the slight additional algorithmic difficulty being more than compensated by the stability, efficiency and flexibility of the Lagrange-Chebyshev interpolant in barycentric form. With this framework in place, we will compare the Lagrange-Chebyshev interpolants to the SBFs and RBFs in Section 4.1.

For all our tests on \textit{open} planar curves, we use nodes contained in $[0,1]$, simply to facilitate the kinds of parametrizations used within the MRS. It is important to find a set of interpolation nodes $\lond_1,\ldots,\lond_{N_d}$ that results in a stable interpolant. When the shape of interest is a closed curve homeomorphic to the circle, evenly-spaced values indeed give us stable interpolants. When moving to an open curve, for stability, we must now choose between different types of clustered nodes mapped to the interval of interest. Regardless of the choice of interpolation nodes, we choose to always \emph{evaluate} our interpolant at a set of equispaced nodes in the interval for static tests. In full Stokeslet simulations involving open curves, we track and interpolate the data sites using a fixed parametric interpolation matrix built out of clustered nodes. The reasoning for this will become apparent when the results of the static tests for different types of nodes are seen. We will explore the effects of different node choices in Section \ref{sec:Results1}.

\section{Using the RBF model within the Regularized Stokeslet method}
\label{sec:RS}

In this section, we will first review the method of Regularized Stokeslets. We will then discuss how both the RBF representations are used (in lieu of finite difference approximations) for geometric modeling within this method.

Consider the scenario where an elastic curve (either open or closed) is immersed in a viscous, incompressible fluid described by the Stokes equations. In two or three dimensions, these equations are as follows:
\begin{align}
\mu \Delta \vu &= \nabla p - \vf,\label{St}\\
\nabla \cdot \vu &= 0,\label{Incomp}
\end{align}
where $\mu$ is the viscosity of the fluid, $p$ is the pressure, $\vu$ is the velocity and $\vf$ is a body force acting on the fluid. For a single point force, the exact solution for the resulting flow is a fundamental solution called a Stokeslet. The Stokeslet can be used to evaluate the fluid velocity at any point off of the structure and is singular if evaluated at the location of the point force. Thus, Cortez derived the method of regularized Stokeslets, removing the singularity by regularizing the forces \cite{Cortez2000,Cortez05}. This method can also handle multiple point forces since the Stokes equations are linear and the solution can be obtained by a superposition of regularized fundamental solutions.

The elastic curve that exerts forces on the surrounding fluid is typically represented by a Lagrangian description, $\bfX(\lambda,t)$ with parameter $\lambda$ at time $t$. 
In the MRS, the continuous formulation of a body force $\vf$ on a closed or open curve $\Gamma$ is written as
\begin{align}
\vf(\vx) = \int_{\Gamma}\vF(\lambda,t) \phi_{\delta}(\vx - \vX)d\lambda,\label{regF}
\end{align}
where $\vF$ is the force per unit length generated along the structure, $\vx$ is any point in the fluid domain (including possibly a point on the structure),  and $\phi_{\delta}$ is a radially symmetric smooth approximation to a delta function. In this formulation, the fluid feels the structure primarily in a region around the curve, governed by the regularization parameter $\delta$. We note that in practice, the Lagrangian structure is discretized into a set of points and the body force $\vf$ is found by approximating the integral in Equation~\eqref{regF} with a  quadrature rule. 

Using a regularized body force as in Equation~\eqref{regF}, the regularized solutions for the pressure and velocity at the IB points (and indeed, in all of space) can be derived. 
If there are $N_s$ forces ${\FF}_k=-\vF_k\bigtriangleup\lambda$  at IB points $\vX_k$ that are equally spaced by $\bigtriangleup\lambda$, the pressure and velocity in $\mathbb{R}^2$ are explicitly given by
\begin{align}
p(\vx) &= \sum_{k=1}^{N_s} \lf(\FF_k \cdot\nabla G_{\delta} (\vx - \vX_k)\rt), \label{eq:pressure}\\
\vu(\vx) &=  \frac{1}{\mu}\sum_{k=1}^{N_s} \lf(\frac{\FF_k}{8\pi}+(\FF_k \cdot \nabla)\nabla B_{\delta}(\vx - \vX_k) -\FF_k G_{\delta}(\vx - \vX_k)\rt). \label{eq:velocity}
\end{align}
$G_{\delta}$ is a smooth approximation of the Green's function and $B_{\delta}$ is a smooth solution to the biharmonic equation such that $\Delta G_{\delta}=\phi_{\delta}$ and $\Delta B_{\delta}=G_{\delta}$. Both $G_{\delta}$ and $B_{\delta}$ are assumed to be radially symmetric, simplifying expressions such as $\nabla B_{\delta}(\vx-\vX)=B_{\delta}'(r)(\vx-\vX)/r$ where $r$ is the Euclidean distance between $\vx$ and $\vX$. In $\mathbb{R}^2$, the blob function  $\phi_{\delta}(\vx-\vX) = \frac{3\delta^3}{2\pi(r^2 + \delta^2)^{5/2}}$ has corresponding regularized functions $G_{\delta}$ and $B_{\delta}'$:
\begin{align}
\label{eq:blob}
G_{\delta}(r) &= \frac{1}{2\pi} \lf[\ln\lf(R + \delta\rt) - \frac{\delta}{R} \rt], \\
B_{\delta}'(r) &= \frac{1}{8\pi} \lf[2r\ln\lf(R + \delta\rt) - r - \frac{2r\delta}{R + \delta} \rt],
\end{align}
where $R=\sqrt{r^2+\delta^2}$.
We remark that taking the limit $\delta \to 0$ for all these expressions will yield the traditional singular expressions for velocity and pressure in the standard Stokeslet formulation. Equations~\eqref{eq:pressure}--\eqref{eq:velocity} are an exact solution for the given regularized forces in Equation~\eqref{regF} and the velocity field in Equation~\eqref{eq:velocity} is everywhere incompressible. For more detail, we refer the reader to~\cite{Cortez2000}.

In the application of ciliary or flagellar movement, an elastic structure is moving and generating forces in a time dependent manner. At each time step, we will have a new body force $\vf$ in the Stokes equations, \eqref{St}-\eqref{Incomp}. As a result of the new force, we can determine the new pressure and velocity field using the MRS and evaluating Equation~\eqref{eq:pressure}--\eqref{eq:velocity}. In the MRS, similar to the IB method \cite{Peskin:2002}, we will assume that the immersed elastic structure is moving at the local fluid velocity. In a standard MRS algorithm, one would evaluate Equations~\eqref{eq:pressure}--\eqref{eq:velocity} at each of the $N_s$ IB points and then update their location by solving a system of $2N_s$ differential equations of the form $\vu(\vX)=\frac{\partial\vX}{\partial t}$.

When using either RBF representation within the method of regularized Stokeslets, the first step is to use the data sites, the appropriate force model and the operators from the previous section to generate forces at the sample sites. We then compute the pressures and velocities at the $N_d$ \emph{data sites} using the forces at the $N_s$ \emph{sample sites}. In other words, we use Equations \eqref{eq:pressure} and \eqref{eq:velocity} but evaluate those equations only at the locations $\vx = \{\vX_j\}_{j=1}^{N_d}$, the set of data sites. While the number of terms ($N_s$) in the summation for the velocity remains the same, we now only need to compute $N_d$ such summations, where $N_d << N_s$. 
By using $N_s$ terms in each of the $N_d$ summations, we maintain the accuracy of the quadrature. The RBF representations thus offer a slight reduction in computational cost, but a great improvement in accuracy for a given computational cost when compared to the traditional algorithm for the method of regularized Stokeslets.

We outline the RBF-Stokeslets algorithm below. 
\begin{itemize}
\item Given a set of data sites $\uvXd$, we use either RBF interpolant to generate a set of sample sites $\uvXs$ (at equispaced parametric nodes) using a variation of Equation~\eqref{eq:nodeeval}.
\item Evaluate forces at the sample sites $\uvXs$ using the differentiation matrices of the RBF interpolant as outlined in Equation~\eqref{diffM}. This will vary based on whether forces depend on tangents, normals, etc., and may involve intermediate steps where derivatives of quantities are computed at data sites.
\item Evaluate the velocity in Equation~\eqref{eq:velocity} at $N_d$ data sites. That is, evaluate $\vu(\vXd_j)$ for $j=1,\ldots,N_d$ where the summation is over the $N_s$ sample sites where $N_d\ll N_s$.
\item Update the location of the $N_d$ data sites by solving a system of $2N_d$ differential equations, $\vu(\vXd_j)=\frac{\partial\vXd_j}{\partial t}$ for $j=1,\ldots,N_d$.
\end{itemize}
After this series of steps, the same process repeats at the next time step. We note again the importance of two sets of sites, data and sampling for the SBF/RBF version of the MRS. The SBF/RBF parametric model of the structure allows efficient and accurate force calculations at $N_s$ sample sites (evenly spaced in the parametric variable $\lambda$). Additionally, the computational expense of summing over point forces in Equation~\eqref{eq:velocity} to update the location of the IB points is decreased by only calculating velocities at $N_d$ data sites. In the results presented in Section 4 for the dynamic elastic filaments, we will use a simple Euler method for the calculation of the updated location of the IB data points. 

\section{Results}

\subsection{Results for static open planar curves}
\label{sec:Results1}

In \cite{Shankar2012}, we tested an SBF interpolant against Fourier-based interpolants on two types of closed objects: a shape represented by an infinitely-smooth function and another shape represented by a twice-differentiable function. In those tests, the SBFs matched the accuracy and convergence of the Fourier interpolants on the smooth object, and were more accurate than Fourier interpolants for the object represented by the twice-differentiable function.

Our goal in this section is to compare the ability of the SBF to approximate twice-differentiable functions against the ability of the standard RBF interpolant to do the same. In addition, we study the behavior of both interpolants against the Lagrange-Chebyshev interpolant. It is important to note that both the RBF and SBF interpolants generalize trivially to 2D parametric interpolation on scattered nodes, while such a generalization of the Lagrange-Chebyshev interpolant is non-trivial.

We start with a simple test function of the form
\begin{align}
\vX_{ideal}(\lon) = (X_{ideal}(\lon),Y_{ideal}(\lon)) = (\lon,b\sin(2\pi\lon - \omega)), \lon \in [0,1].
\label{eq:obj_ideal}
\end{align}
This test function represents an idealized sinusoidal shape, and is often the parametrization used for a sperm flagellum within regularized Stokeslet simulations \cite{Gillies09,Olson11,Olson13b,Simons14}. Obviously, if this shape always stayed sinusoidal, it would be impossible to have a better interpolant than a parametric Fourier interpolant, which contains a sine term in its expansion and can consequently reproduce the sinusoidal wave with spectral accuracy. However, over the course of a simulation, it is likely that the immersed elastic curve may deviate from the ideal sinusoidal shape. The resulting shape may be represented either by a smooth function or by a finitely-differentiable function, depending on the nature of the spatial discretization of the regularized Stokeslet method. We now consider the case of twice-differentiable functions.

\begin{figure}[ht]
\centering
\includegraphics[height=2.5in,width=2.7in]{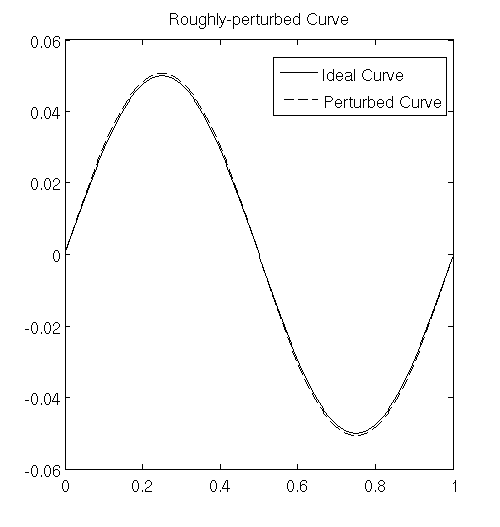}
\caption{Test object used for error analysis. The ideal shape is shown with a solid line (-) and the roughly perturbed version is shown with a dashed line (- -). Note the visual similarity of the perturbed shape to the ideal shape despite the difference in their mathematical representations.}
\label{fig:shape2d}
\end{figure}

We define a perturbation of this idealized shape (see Figure \ref{fig:shape2d}). This perturbation results in a curve that, while visually not very different from the idealized shape, has only two underlying derivatives. We only perturb the y-components of the shape for this test. The shape is given by:
\begin{align}
Y_{P} = \lf[1.0 + A \exp\lf(\frac{\lf( - \lf( 1-\cos^{2}(2\pi\lon) \rt)^{1.5} \rt)}{\sigma} \rt)\rt] \, Y_{ideal} \label{eq:obj_2d}.
\end{align}
In the above equation, $A$ controls the amplitude of the perturbation, with $\sigma$ controlling the width of the perturbation. For the tests that follow, we use $A = 0.04$ and $\sigma = 0.9$. This allows us to reproduce the conditions of the test in \cite{Shankar2012}, now allowing us to compare alternatives to Fourier interpolants. For all our tests, we set $b = 0.05$ in Equation \eqref{eq:obj_ideal}. These are the values used for the plot in Figure \ref{fig:shape2d}. 

Since we are parametrically interpolating open curves, it is unwise to use equispaced nodes lest we encounter the Runge phenomenon. We therefore present results only for clustered nodes. However, since Fourier interpolants are likely unstable on clustered nodes, and are inappropriate for the modeling of curves represented by twice-differentiable functions, we do not present any results for this case. In our experiments, the poor accuracy of Fourier interpolants in interpolating such shapes was quite similar to that seen in \cite{Shankar2012}.

Before we proceed, we quickly remark on the choice of interpolant for infinitely-smooth open curves. For curves that are represented by periodic and smooth functions, the obvious victor is the Fourier interpolant, just as in the closed curve case seen in~\cite{Shankar2012}. However, similar to the results from that work, SBF interpolants are competitive with the Fourier interpolants in this scenario, with the only drawback being that errors cannot be pushed all the way down to machine precision for the SBF interpolants due to ill-conditioning. This may be somewhat surprising, since the use of SBFs implies that we are parametrizing open curves on the circle. However, we have found that this approach works perfectly well. Of course, since the curves are open, SBFs have no advantage over standard RBFs (unlike in the closed case, where RBFs fail to capture the periodic point well). However, SBFs appear to be equivalent to RBFs in this setting. Depending on the performance of SBF interpolants in the following sections, this implies that one could write very general software based on SBF interpolants that would work excellently for the modeling of both closed and open curves. We will not present results for the recovery of infinitely-smooth functions corresponding to open curves with either SBFs or RBFs in this work, and focus on the more interesting case of twice-differentiable functions representing open curves.

For the modeling of open, twice-differentiable curves with RBF and SBF interpolants, we use two types of clustered parametric nodes: standard Chebyshev nodes, and the so-called ``mapped Chebyshev nodes'' introduced by Kosloff and Tal-Ezer. For brevity, we will refer to the mapped Chebyshev nodes as KTE nodes. In \cite{PlatteAnalytic}, it was shown that the Lebesque constants corresponding to RBF interpolation at the KTE nodes grow much more slowly than in the case of RBF interpolation at Chebyshev nodes. For a parametric Chebyshev node $\lon^{cheb}$ defined in Equation \eqref{eq:cheb_nodes}, the corresponding KTE node is given by 
\begin{equation}
\lon^{KTE} = \frac{\sin ^{-1}(\alpha \lon^{cheb})}{\sin ^{-1} (\alpha)}.
\end{equation}
When $\alpha = 1$, the nodes are equispaced, while the limit $\alpha \to 0$ recovers the Chebyshev nodes. 

\subsubsection{Interpolation at Chebyshev nodes}

\begin{figure}[htbp]
\centering
\includegraphics[scale=0.5]{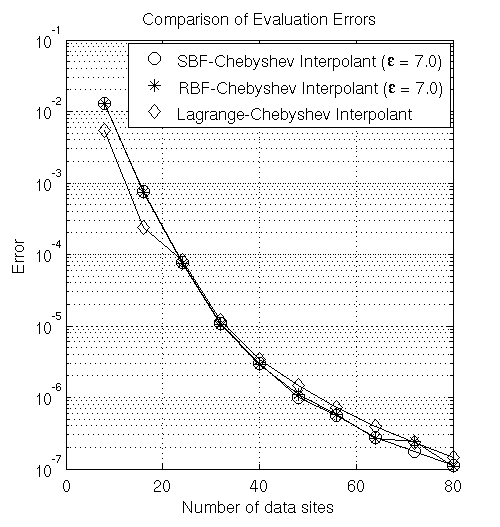}
\includegraphics[scale=0.5]{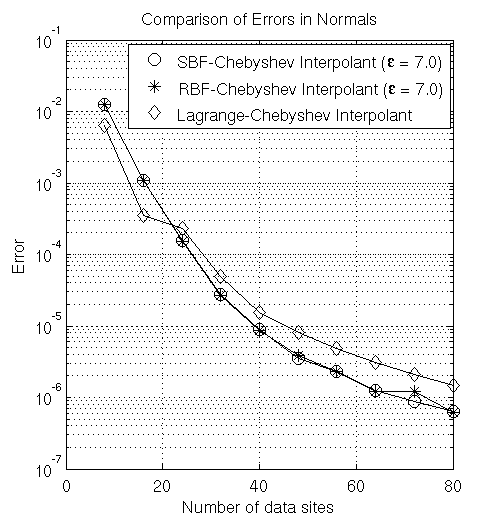}
\includegraphics[scale=0.5]{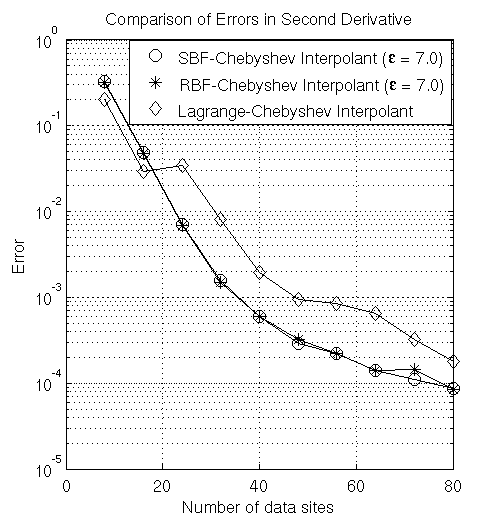}
\caption{Comparison of interpolants at Chebyshev nodes when interpolating and evaluating a shape represented by a twice-differentiable function (top left), its unit normals (top right) and its second derivative (bottom). The shape is given by Equation \eqref{eq:obj_2d}.}
\label{fig:cheb_comp}
\end{figure}

We first compare the ability of the SBF, RBF and barycentric Lagrange interpolants to approximate $Y_P$, its normals and its second derivatives with respect to $\lon$. All interpolants are built on the same set of Chebyshev data sites nodes in $(0,1)$, and all interpolants are evaluated at a set of $N_s = 400$ equispaced sample site nodes in $[0,1]$. The results are shown in Figure \ref{fig:cheb_comp}. For each case in Section 4.1, the error is defined as the maximum of the point-wise $\mathcal{L}^2$ errors for each of  the sample sites.

From Figure \ref{fig:cheb_comp}, we see that both the SBF and the RBF have very similar errors in all cases for the same shape parameter at the Chebyshev nodes. The figure on the top left shows that both the SBF and RBF interpolants have errors that are quite similar to the errors of the barycentric Lagrange interpolant, albeit slightly lower for higher values of $N_d$. However, the figure on the top right and the one on the bottom show that the SBF and RBF interpolants have lower errors when computing normals, and much lower errors for the second derivatives. This may be a consequence of the difficulty in interpolating a function and evaluating derivatives at a different set of nodes in the barycentric formulation of the Lagrange-Chebyshev interpolant. Interestingly, the SBF seems to have very similar errors to the standard RBF, despite the difference in the distance argument passed in to the basis function $\phi(r)$. While SBFs are preferred for the interpolation of analytic and/or very smooth periodic functions corresponding to closed curves (as seen in \cite{Shankar2012}), this indicates that they can safely be used in the context of modeling open, twice-differentiable shapes with high accuracy as well.

\subsubsection{Interpolation at KTE nodes}

\begin{figure}[htbp]
\centering
\includegraphics[scale=0.5]{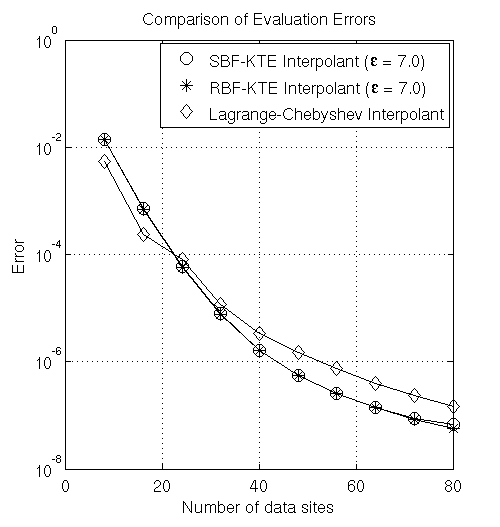}
\includegraphics[scale=0.5]{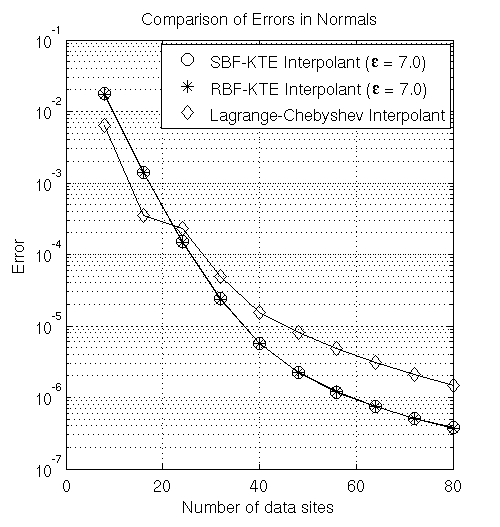}
\includegraphics[scale=0.5]{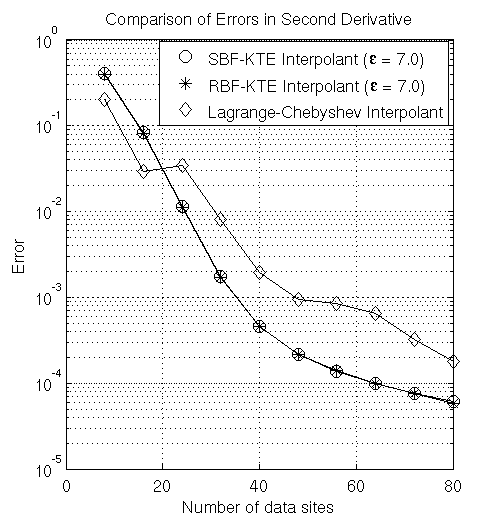}
\caption{Comparison of SBF/RBF interpolants at KTE nodes, with the Lagrange interpolant at Chebyshev nodes, when interpolating and evaluating a shape represented by a twice-differentiable function (top left), its unit normals (top right) and its second derivative (bottom). The shape is given by Equation \eqref{eq:obj_2d}.}
\label{fig:kte_comp}
\end{figure}

We now compare the SBF and RBF interpolants built on KTE nodes against the barycentric Lagrange-Chebyshev interpolant in the approximation of $Y_P$, its normals and its second derivatives. The KTE nodes were generated with $\alpha = 0.85$ and a shape parameter of $\ep = 7.0$. Note that at this value of $\alpha$, one cannot safely use the Lagrange interpolant with the KTE nodes, since these nodes are quite close to evenly-spaced.  The results are shown in Figure \ref{fig:kte_comp}.

A careful examination of Figure \ref{fig:kte_comp} shows that both the SBF and RBF interpolants have even lower errors at KTE nodes than they do at Chebyshev nodes. Consequently, they have much lower errors than the barycentric Lagrange-Chebyshev interpolant for evaluation, normal computation and computation of second derivatives. Once again, the SBF and RBF interpolants have almost identical errors for a given value of $N_d$ when interpolating this twice-differentiable function. This shows the promise of the SBF in interpolating functions of very different smoothness.

We note that even for infinitely-smooth target functions (results not shown), the SBF-KTE interpolant can beat both the SBF interpolant at equispaced nodes and the RBF interpolant at either set of nodes. For some finite and small value of $N_d$, the SBF-KTE interpolant actually has lower errors than the Fourier interpolant (the latter collocated at an equispaced set of nodes). The errors for the Fourier interpolant, unlike those for the SBF-KTE interpolant, can be pushed all the way down to machine precision. To do so for the SBF interpolants, we would require a so-called ``flat'' algorithm (these tend to be more expensive), such as the RBF-QR algorithm~\cite{FornbergPiret:2007}.

Interestingly, in our experiments with the RBF-QR algorithm, we have observed that as $\ep \to 0$, the SBF-KTE interpolant becomes unstable, and a change to the standard equispaced nodes is required; similarly, the RBF-KTE interpolant must be swapped out for the RBF-Chebyshev interpolant as $\ep \to 0$. Both these changes are to be expected, since in that limit, the SBF and RBF interpolants recover the Fourier and Lagrange interpolants respectively. We do not explore such strategies here, as our goal is to attain a cost-effective and highly accurate approximation of both function and derivative approximation for small values of $N_d$ and non-zero values of $\ep$.

\subsubsection{Interpolation at KTE nodes with ``best'' shape parameters}

\begin{figure}[htbp]
\centering
\includegraphics[scale=0.45]{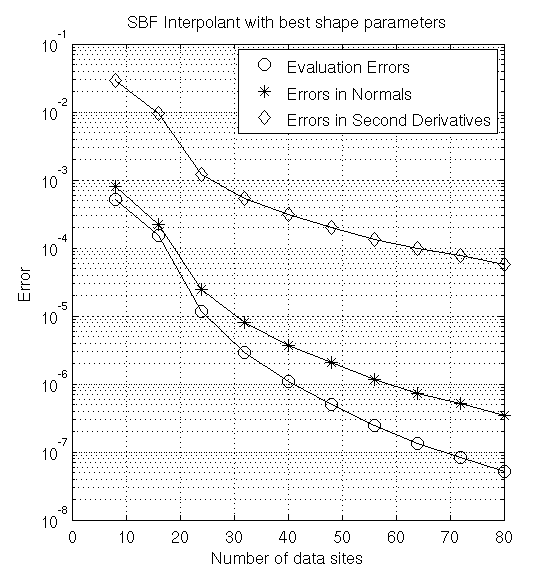}
\includegraphics[scale=0.45]{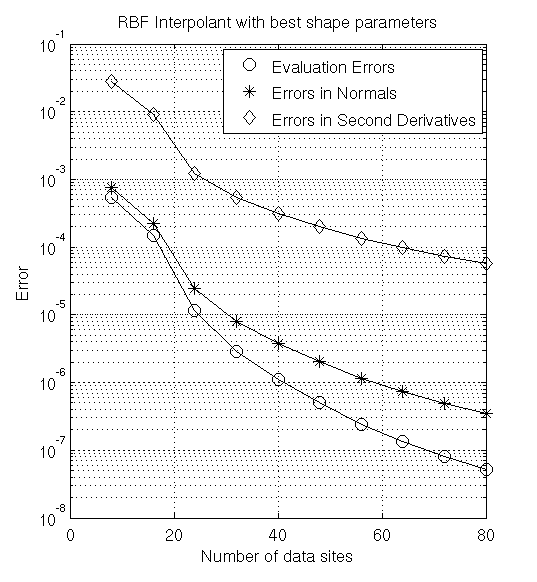}
\caption{Errors in evaluation, normal computation and calculation of second derivatives for the SBF (left) and RBF (right) interpolants for carefully selected shape parameters.}
\label{fig:best_comp}
\end{figure}

Until now, we have made the rather simple choice of selecting a single shape parameter for all values of $N_d$ for both the SBF and RBF interpolants. This choice will negatively affect the errors for smaller values of $N_d$. Indeed, for a given value of $N_d$, the correct strategy is to pick a single shape parameter that results in the lowest error. It is possible to do so for each value of $N_d$ in an error plot. We run such an experiment at the KTE nodes (which we have already determined to be a better node set than the Chebyshev nodes for the SBF and RBF) using $\alpha=0.85$. The results of this test are shown in Figure \ref{fig:best_comp}. 

As expected, both the SBF and RBF interpolants now have lower errors for the smaller values of $N_d$ than they did when using the fixed value of $\ep = 7.0$. Furthermore, while the error curves do not look quite as pleasing as those in Figure \ref{fig:kte_comp}, it is clear that even with these optimized shape parameters, the SBF and RBF interpolants converge reasonably well. In addition, the errors for the SBF and RBF interpolants are now lower than those of the Lagrange-Chebyshev interpolant even for small values of $N_d$, which was not the case for the studies with fixed $\ep$. While there exist a few sophisticated algorithms to optimize the shape parameter, we simply pick the best ones by sampling densely in $\ep$-space, and computing errors for all those values of $\ep$ for a given number of data sites $N_d$. For reproducibility, we document the results in Table \ref{tab:eps}.

\begin{table}[htbp]
\centering
\begin{tabular}{|c|c|c|}
\hline
$N_d$ & $\ep_{SBF}$ &  $\ep_{RBF}$ \\ \hline
8 &  2.5 & 2.6  \\ \hline
16 &  3.2 & 3.2 \\ \hline
24 &  3.0 & 2.9 \\ \hline
32 &  3.8 &  3.9 \\ \hline
40 &  3.8 &  3.6 \\ \hline
48 &  4.5 &  4.7 \\ \hline
56 &  5.7 &  5.9 \\ \hline
64 &  8.2 &  8.0 \\ \hline
72 &  8.9 &  8.8 \\ \hline
80 &  9.0 &  8.6 \\ \hline
\end{tabular}
\caption{``Best'' values of the shape parameter $\ep$ corresponding to different numbers of data sites $N_d$ for the SBF and RBF interpolants. $100$ samples of $\ep$ were used as the candidate set.}
\label{tab:eps}
\end{table}

While the values in Table \ref{tab:eps} look highly tuned, this is simply an artifact of the range of samples used to generate the different values of $\ep$. For a given $N_d$, one can certainly use values of the shape parameter that are close to but different from those presented in the table, with little change in the error. This table is merely intended as a guide to selecting the shape parameter in preparation for a full time-dependent simulation, and for demonstrating the value in selecting smaller shape parameters for smaller $N_d$.

\subsection{Stokeslet simulations on closed shapes}
\label{sec:Results3}

In this section, we show results for simulations of closed elastic objects immersed in a viscous, initially stationary fluid. Error analysis of RBFs with closed curves has been previously detailed in the context of the IB method \cite{Shankar2012}, thus our results focus on showing that this method works for closed curves in the context of the regularized Stokeslet method. First, we present results for a stationary problem that has an exact solution. We then show results for a time-dependent fluid-structure interaction simulation. For both test cases, we use the multiquadric (MQ) radial kernel function given in Equation \eqref{eq:mq} and the SBF given in Equation \eqref{SBFr} to interpolate the y-coordinate at the data sites. 

\subsubsection{Computation of tangential forces}\label{sec:Cortez}
\begin{figure}[ht]
\centering
\includegraphics[scale=0.4]{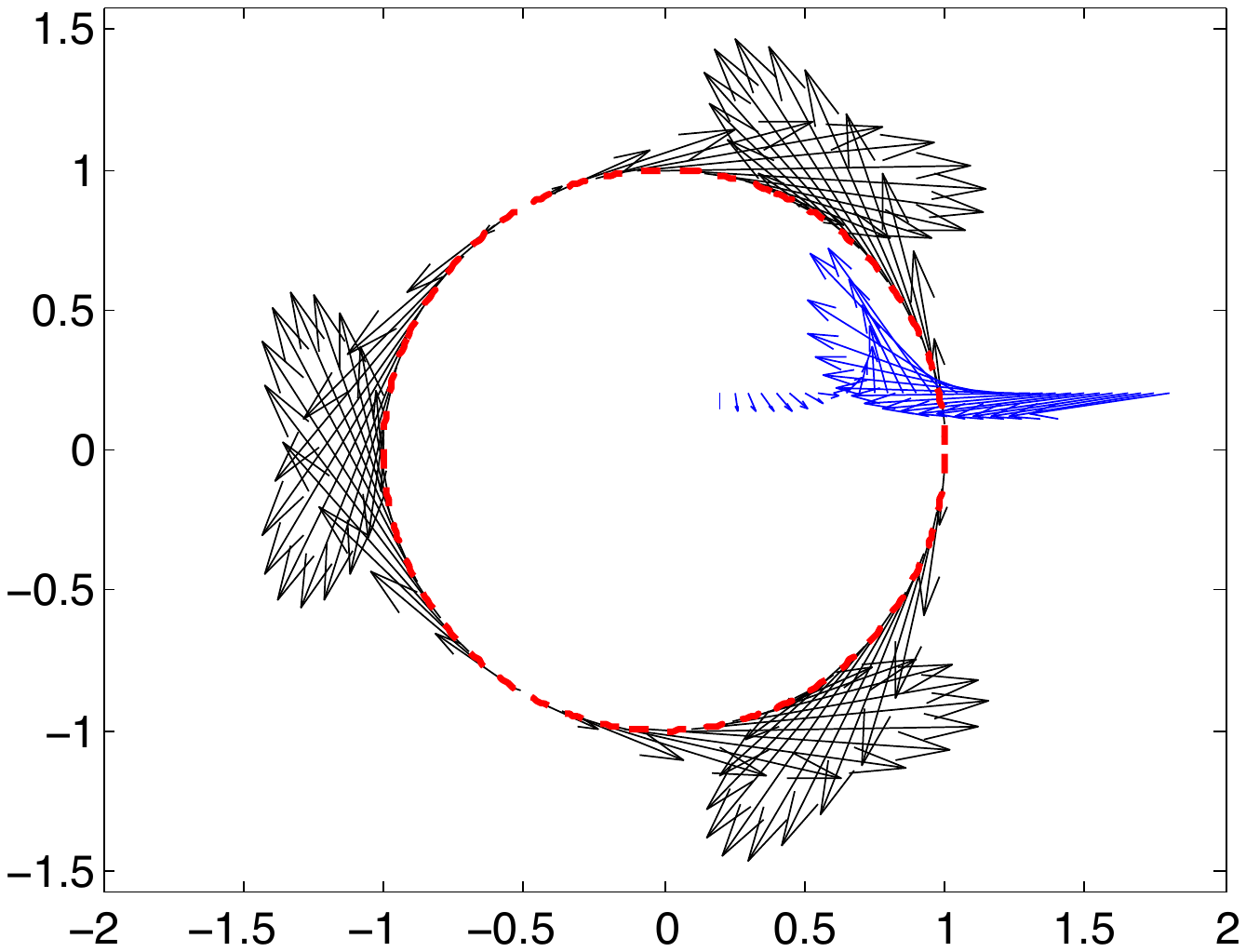}
\includegraphics[scale=0.4]{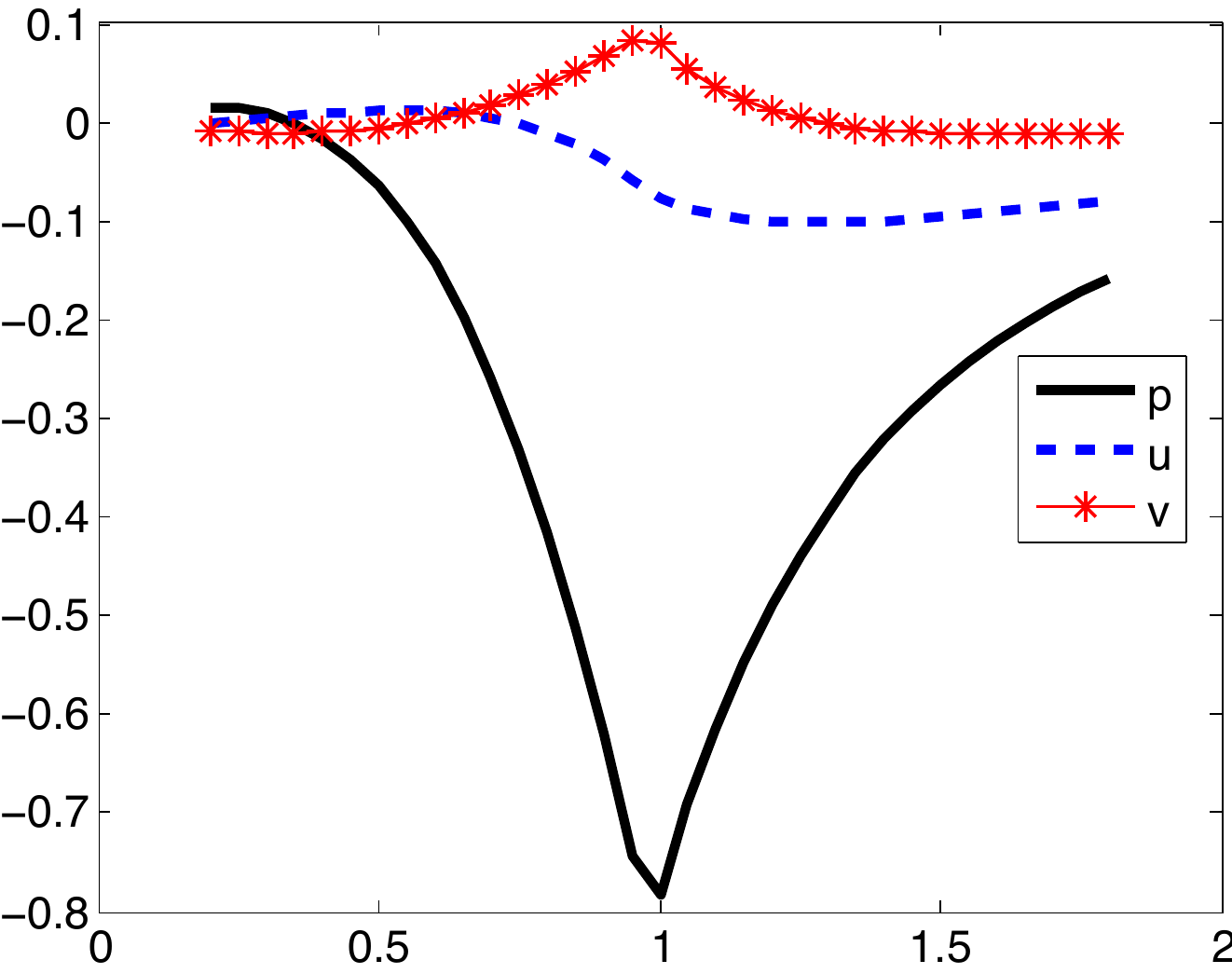}
\caption{Exact solution for a unit circle with tangential forces, presented in~\cite{Cortez2000}. The left figure shows the structure, the unit circle, with forces on the structure shown with black arrows. The line of arrows correspond to evaluating the fluid velocity at the points on the line ($x$,1/5). On the right, the exact solution for the pressure $p$ and velocity $\vu=(u,v)$ evaluated at $y=1/5$ for the $x$ values displayed. }
\label{fig:ExactCortez}
\end{figure}

This test problem was presented in~\cite{Cortez2000} (example $4b$). The forces $\vF$ on the an elastic structure corresponding to the unit circle are prescribed to be purely tangential to the boundary,
\begin{equation}
\vF(\lon)=-2\sin(3\lon)\frac{\partial \vX(\lon)}{\partial \lambda} \label{Fexact}
\end{equation}
for parametric variable $\lambda$. The unit circle and corresponding tangential forces (black arrows on the structure) are shown on the left plot in Figure \ref{fig:ExactCortez}. The exact solution for Stokes equation with viscosity $\mu=1$, as detailed in~\cite{Cortez2000} and shown on the right in Figure \ref{fig:ExactCortez}, has a discontinuous pressure and velocity gradient ($\nabla p$ and $\nabla \vu$) across the boundary of the elastic structure and continuous $p$ and $\vu$ everywhere in the infinite 2D fluid. 

\begin{figure}[ht]
\centering
\includegraphics[scale=0.4]{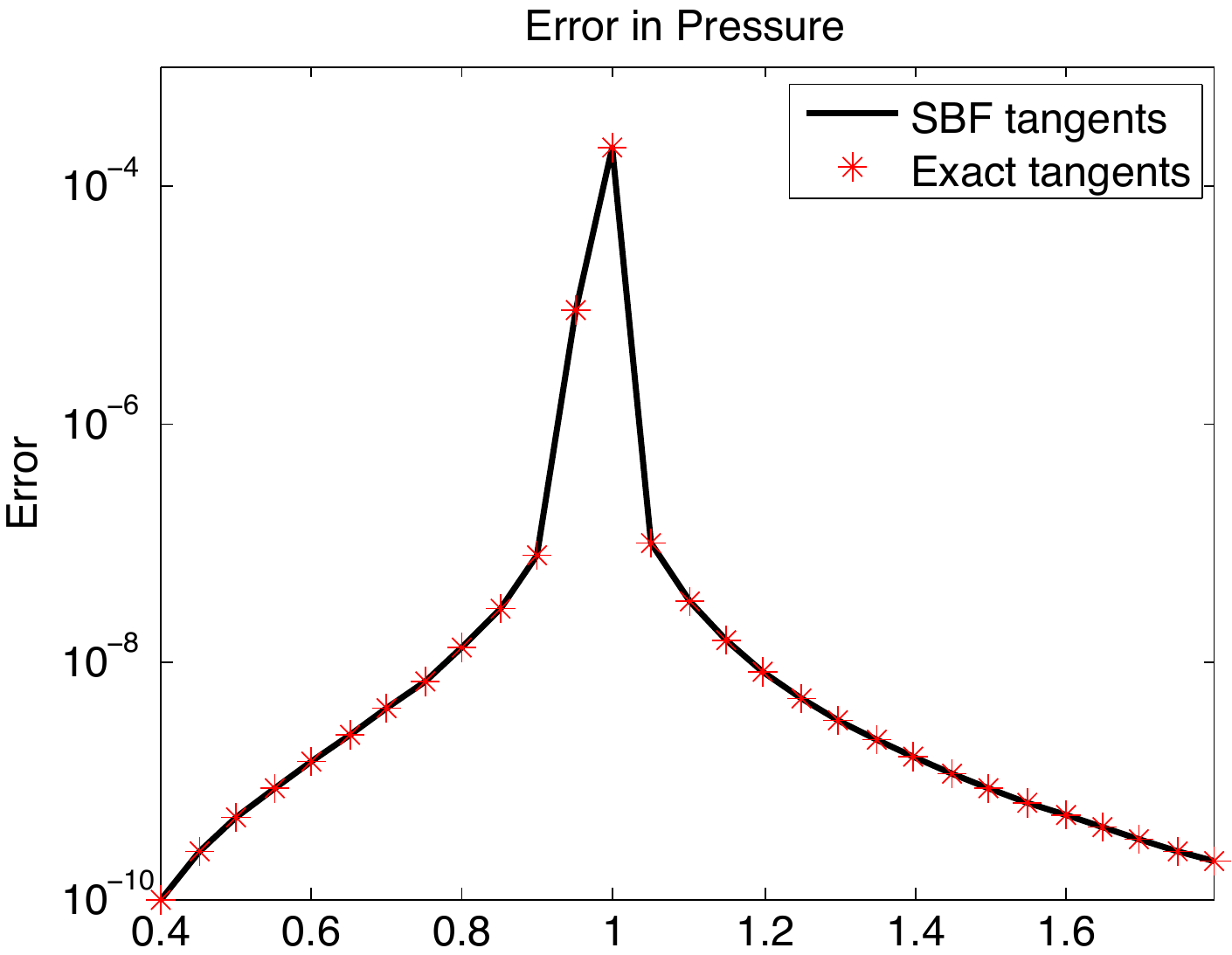}
\includegraphics[scale=0.4]{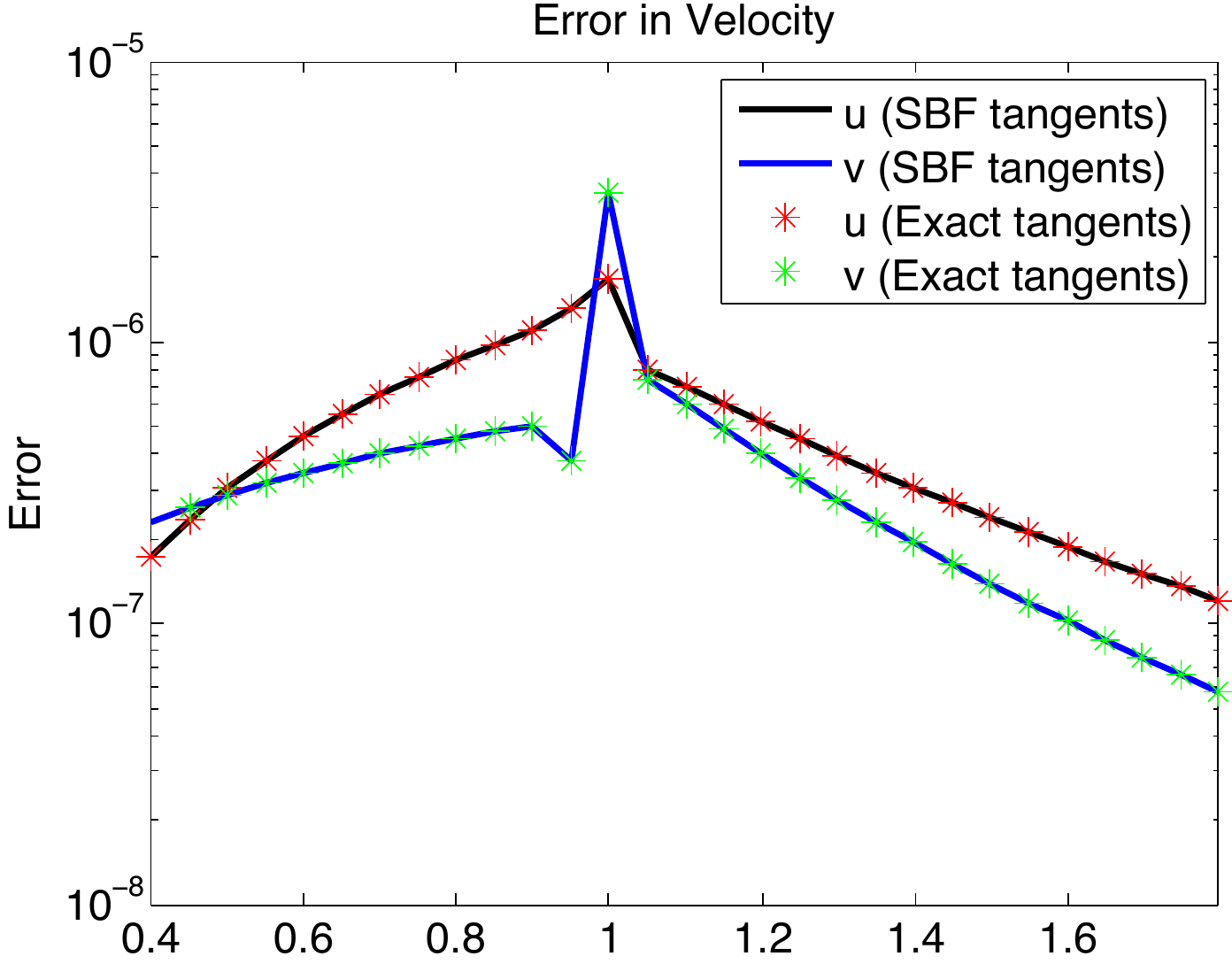}
\caption{Errors in the pressure and velocity fields when using the SBF model. The figure on the left shows the errors in the pressure for the test case presented in~\cite{Cortez2000}, while the figure on the right show the errors corresponding to each component of the velocity field. $N_d = 25$ data sites and $N_s = 400$ sample sites are used for the SBF representation and 400 sites are used for the exact force evaluation with the MRS. The numerical solutions are compared against the exact solution given in~\cite{Cortez2000}.}
\label{fig:CortezTest1}
\end{figure}

We use $N_d = 25$ equispaced data sites for the representation of the circle with $\lon\in[0,2\pi)$. The derivative matrix in Equation \eqref{diffM} is used to compute the tangents ${\partial \mathbf{X}(\lambda)}/{\partial \lambda}$ on the circle in order to calculate the force in Equation~\eqref{Fexact}. We evaluate the forces at $N_s = 400$ equispaced sample sites, and then compute the pressure and velocity using Equations \eqref{eq:pressure}--\eqref{eq:velocity} at a set of markers lying on the line $(x,\frac{1}{5}),0.4\leq x\leq 1.8$. We compare the errors in the solution obtained from the SBF representation with those in the solution obtained by using the analytically correct tangents on the circle. The error reported is the absolute value of the difference of the solutions obtained by using the two methods. The results are shown in Figure \ref{fig:CortezTest1}.

The results in Figure \ref{fig:CortezTest1} show that despite using only $N_d = 25$ data sites to represent the structure with SBFs and sampling the forces at $N_s = 400$ sample sites, the absolute value of the error in the SBF-based Stokeslet approximations to the pressure and velocity fields are almost identical to those obtained when using analytically computed tangents at $N_s = 400$ IB points using MRS. For both methods, the error decreases as you move away from the boundary at $x=1$. Similar results were also observed for the RBF-based Stokeslet approximations. On the circle, one can expect similar performance from a Fourier interpolant, though the Fourier interpolant will offer lower accuracy on roughly-perturbed closed shapes for the same number of data sites~\cite{Shankar2012}.

\begin{figure}[ht]
\centering
\includegraphics[scale=0.4]{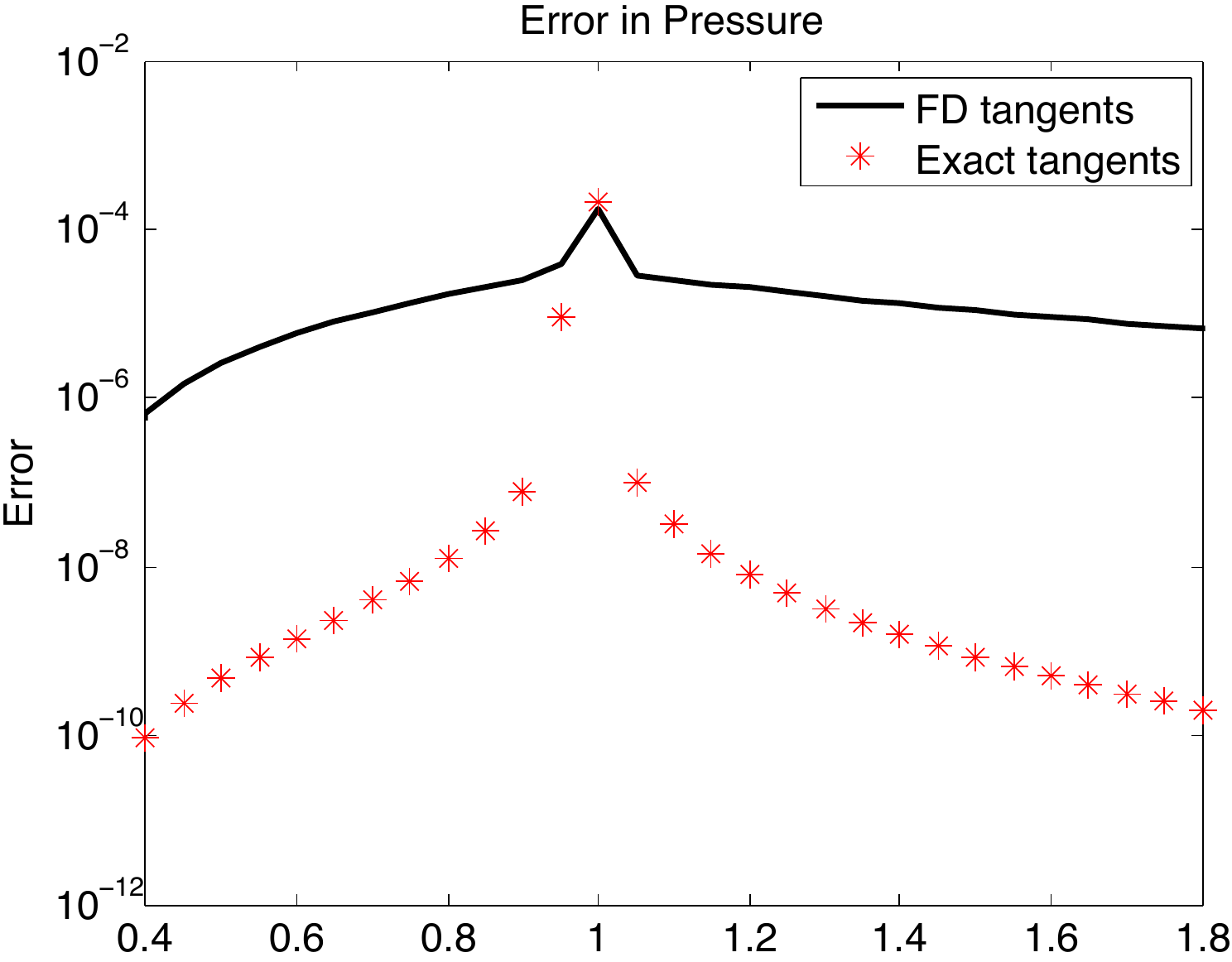}
\includegraphics[scale=0.4]{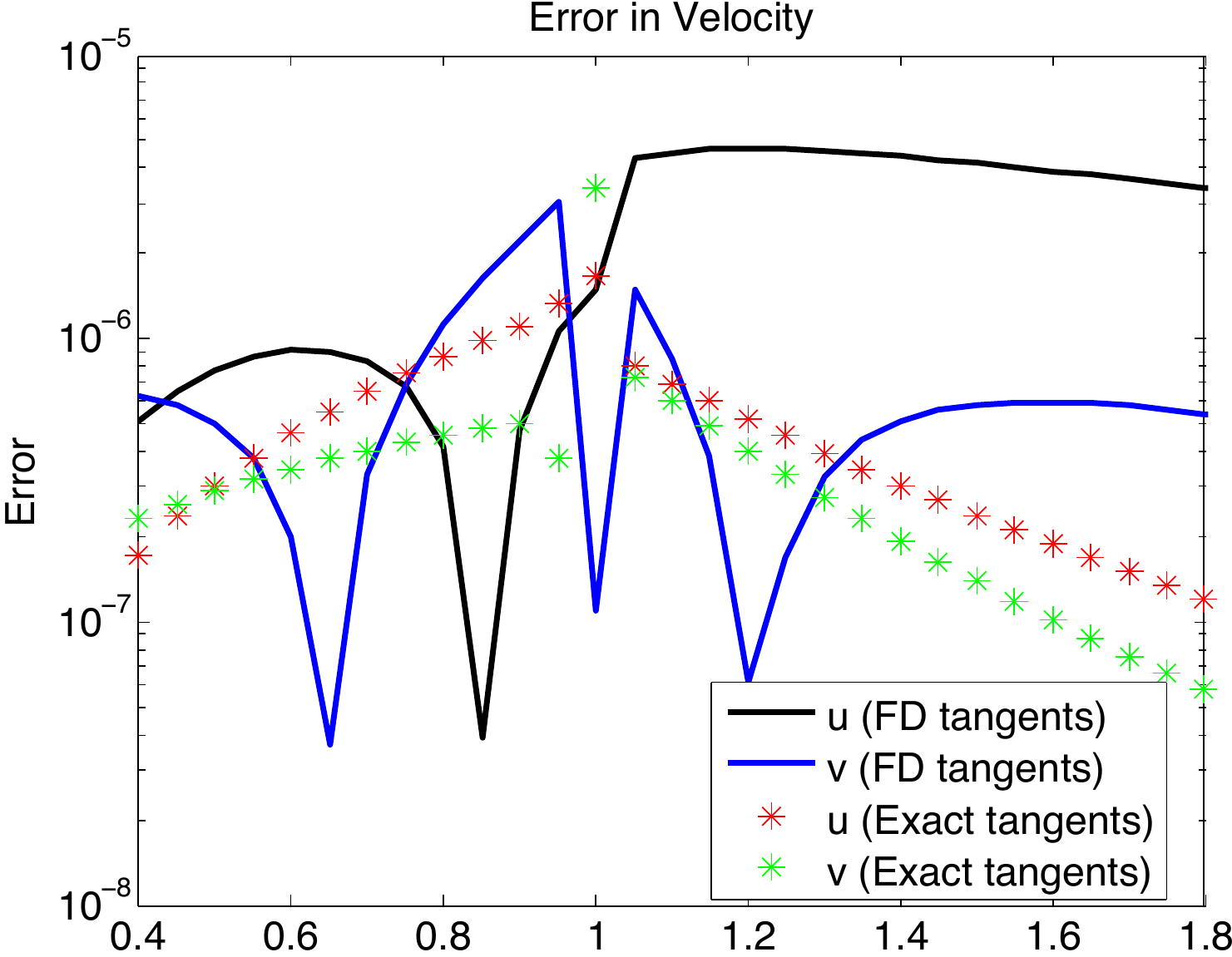}
\caption{Errors in the pressure and velocity fields when using finite difference approximations to the tangents. The figure on the left shows the errors in the pressure for the test case presented in~\cite{Cortez2000}, while the figure on the right show the errors corresponding to each component of the velocity field. $N_s = 800$ IB points are used for the finite differences. The numerical solutions are compared against the exact solution given in~\cite{Cortez2000}.}
\label{fig:CortezTest2}
\end{figure}

To understand the power of global representations, we repeat the test with a second-order finite difference approximation to the tangent. In the Stokeslet and IB literature, finite difference approximations to derivatives at IB points are often computed using differences between points on the boundary that are halfway between the IB points (\emph{e.g.},\cite{OlsonLimCortezJCP2013}). This is equivalent to using $N_s = 800$ IB points (or $N_s = 400$ IB points with ``tight'' finite difference stencils). Figure \ref{fig:CortezTest2} shows the results for the pressure and velocity obtained with the finite difference approximation at $N_s = 800$ IB points. Clearly, the error for the pressure is significantly higher when using the finite difference approximations to the tangents, in contrast to the results achieved using $N_d = 25$ data sites and $N_s = 400$ sample sites with the SBF representation. On the right, the finite difference approximations at $N_s=800$ leads to fairly low error, but it oscillates dramatically and has mean square error that is greater than using the exact formula for the tangents. This highlights a hazard of using finite differences to approximate geometric quantities. In comparison, the SBF solution is smooth and oscillation-free. Additionally, we emphasize that we are able to achieve a lower error in the SBF-stokeslets method than the FD method when using only $N_d=25$ data sites and $N_s=400$ sample sites. We remark that~\cite{Cortez2000} employs correction terms to the method of regularized Stokeslets on this test case to improve its accuracy. We do not use that approach in this current work.

\subsubsection{Time-dependent simulations.}\label{closedSim}
We will use a closed curve as our first test case, as previously used in other methods \cite{Cortez00,Cortez12}. We will start with a closed, elastic structure that is initialized as a perturbation of a circle and immersed in the fluid. The initial configuration $\bfX(\lambda,0)$ is defined as
\begin{equation}
\vX(\lambda,0)=\left((1+\beta\cos(\nu\lambda))\cos(\lambda),(1+\beta\cos(\nu\lambda))\sin(\lambda)\right)
\end{equation}
for $\lambda\in [0,2\pi)$ with $\nu$ and $\beta$ as parameters. For this closed case, we will use $N_d=25$ data sites and $N_s=50$ sample sites; each of these will be equally spaced points in the parameter $\lambda$. The force will be proportional to the curvature at a given time point and defined as follows,
\begin{equation}
\vF(\lambda,t)=-\frac{1}{10}\kappa(\lambda,t)\unrml(\lambda,t)\left(\mathcal{L}(t)-\frac{3\pi}{2}\right),
\label{closedF}
\end{equation}
where $\kappa(\lambda,t)$ is a signed curvature, $\unrml$ is the unit normal, and $\mathcal{L}(t)$ is the arclength of the structure at time $t$. These forces will tend to restore the shape of the closed elastic structure to that of a circle.  

\begin{figure}[ht]
\centering
\hspace{0.4cm}\includegraphics[scale=0.4]{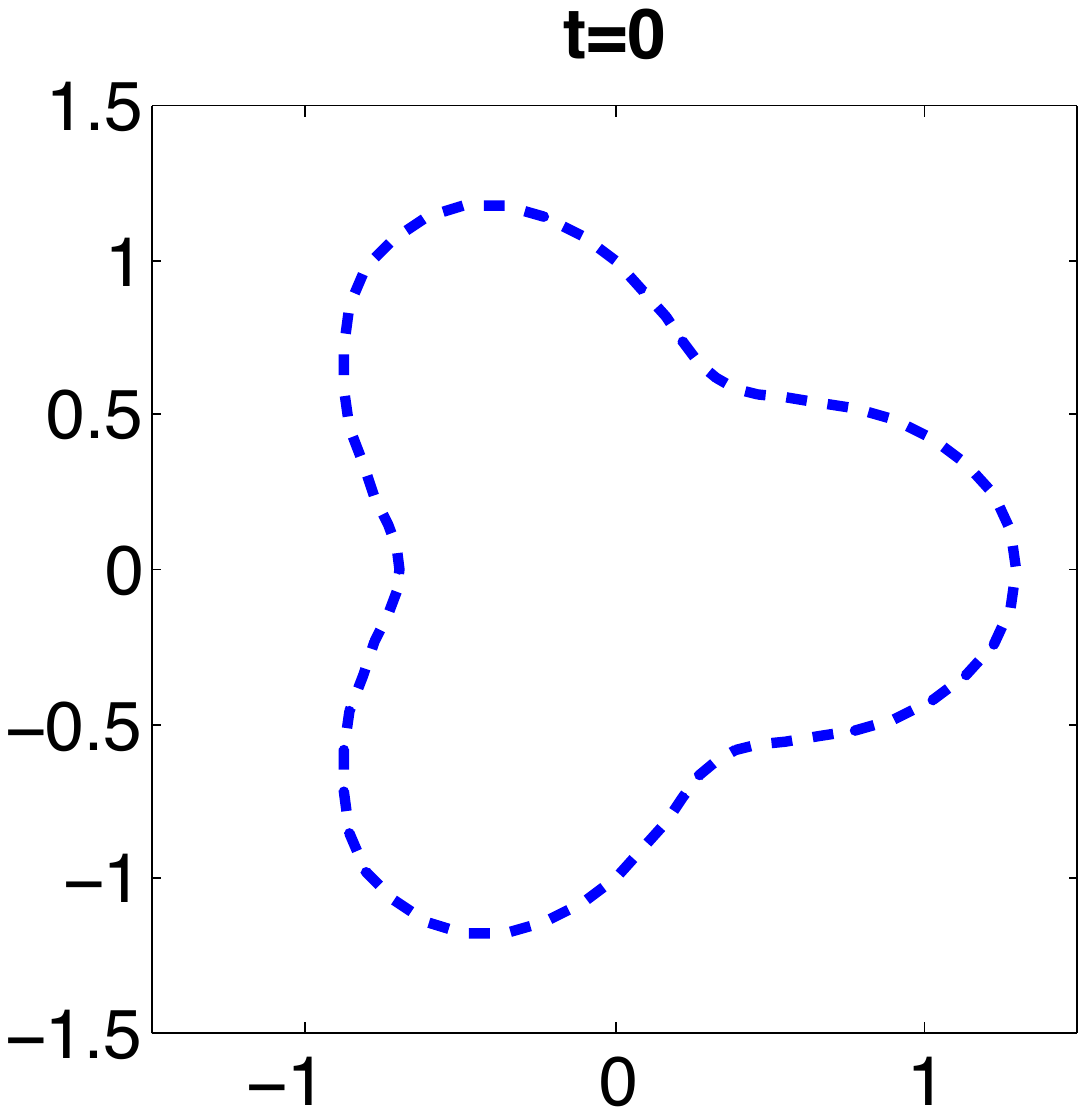}
\hspace{1.0cm}\includegraphics[scale=0.4]{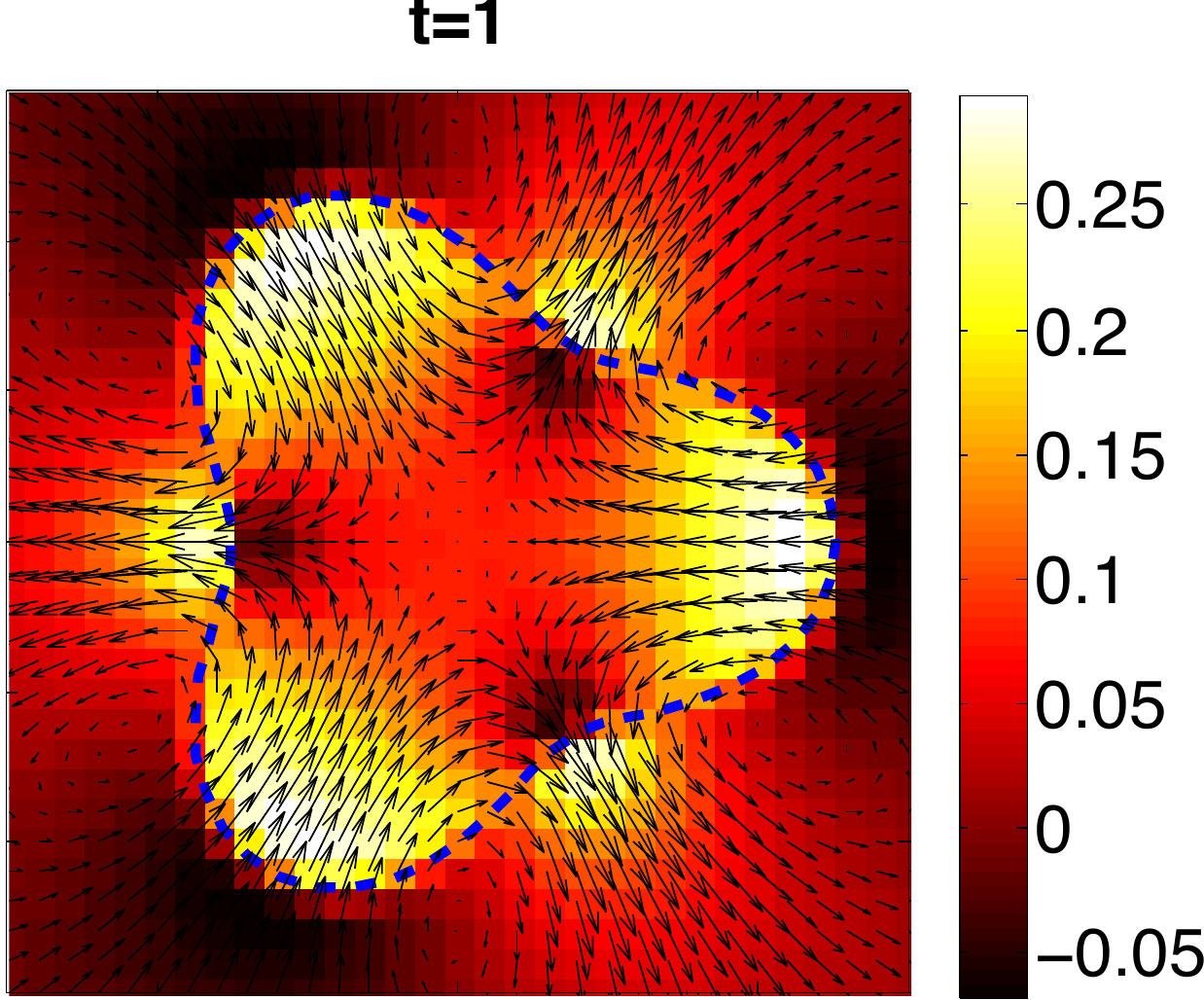}\\
\includegraphics[scale=0.4]{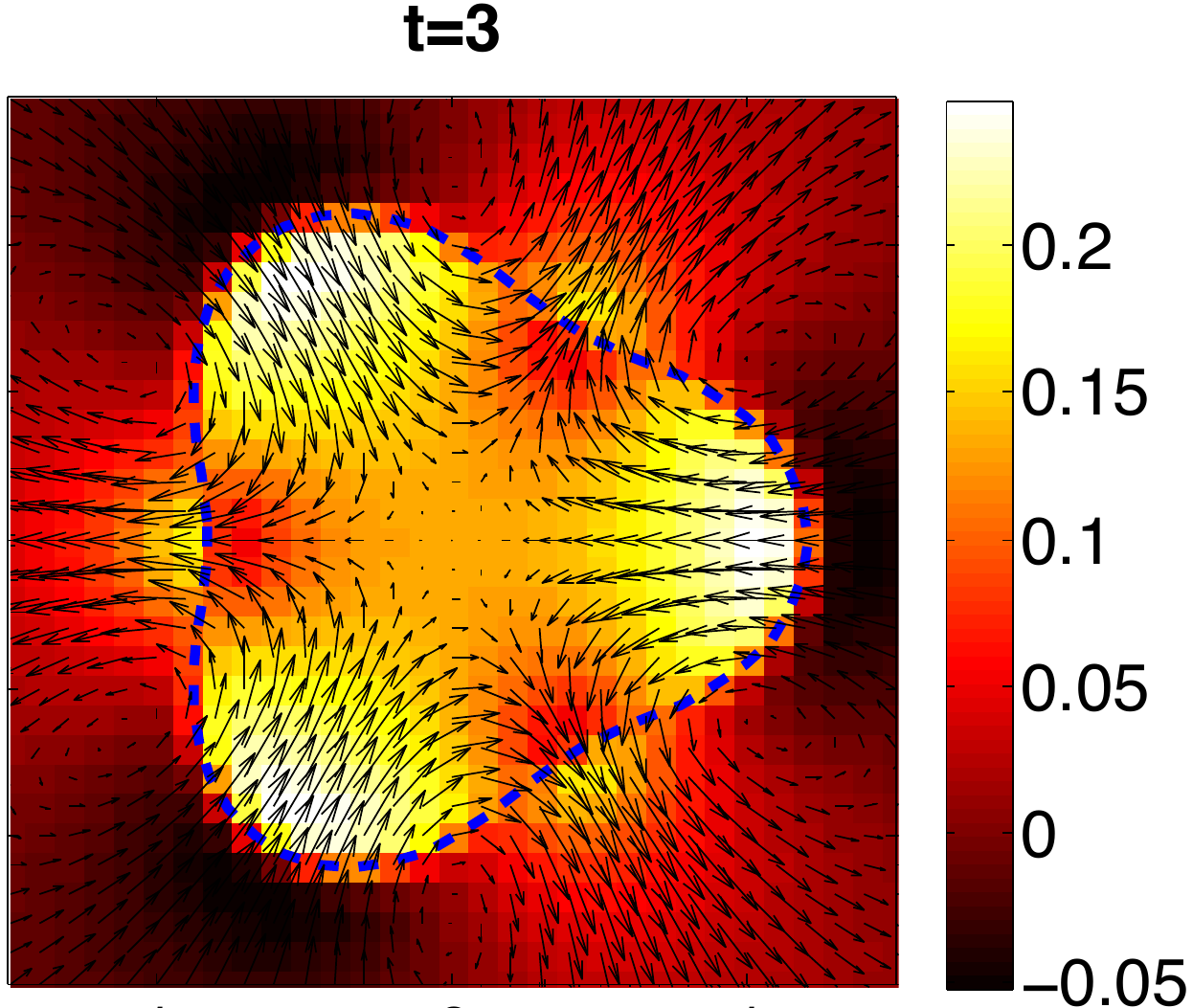}
\includegraphics[scale=0.4]{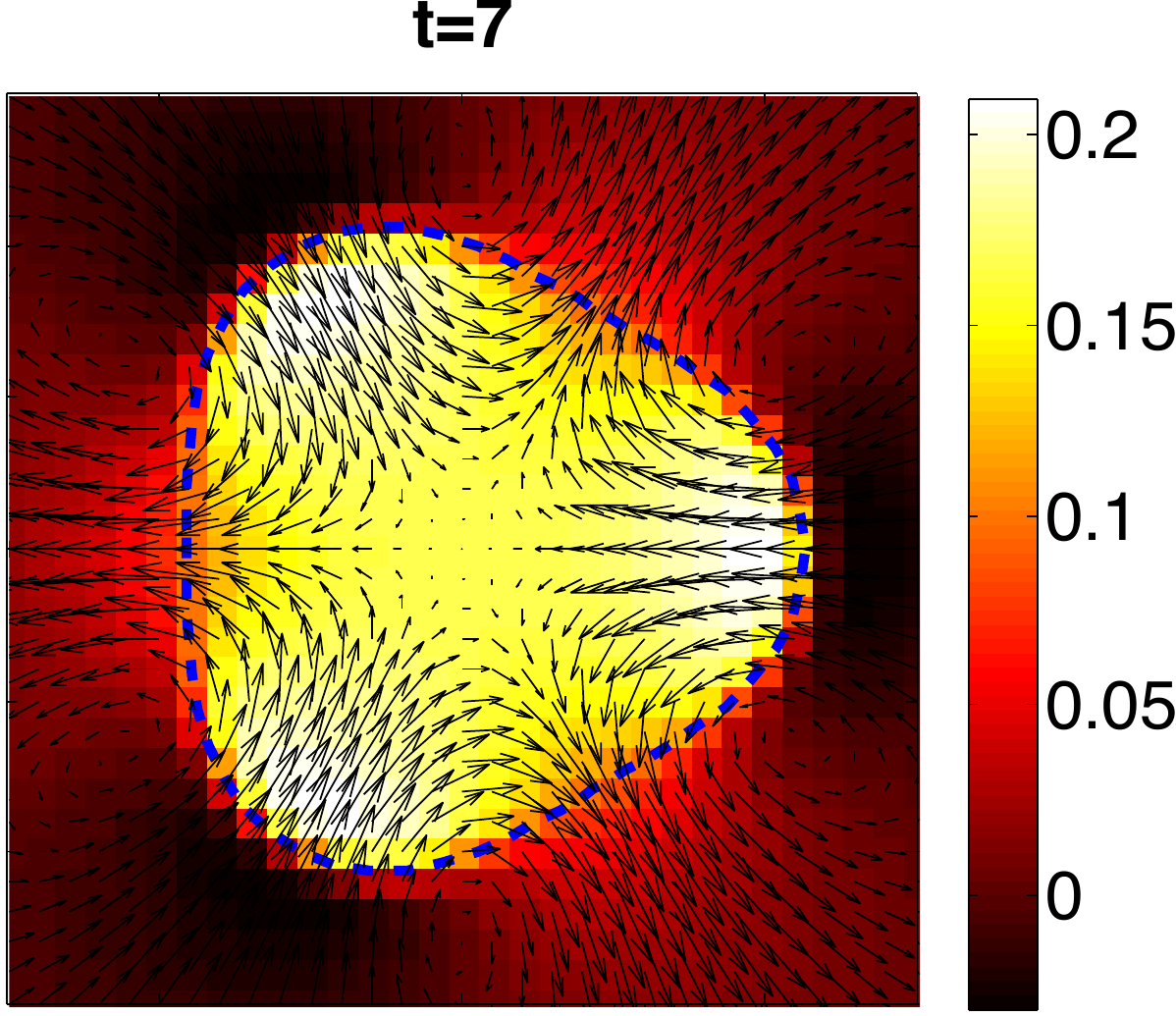}\\
\caption{Results for a closed elastic structure initialized with $\beta=0.3$, $\nu=3$, $N_d=25$, $N_s=50$, $\Delta t=10^{-3}$, $\delta=4\pi/N_s$, and $\ep = 1.1$. Upper left is the initial configuration and the other three panels show the flow field with the black arrows and the color bar corresponds to the pressure. }
\label{fig:RSBTtest}
\end{figure}

To compute the time dependent movement of the elastic structure, we will use Equation \eqref{diffM} to evaluate the derivatives of the SBF y-coordinate interpolant at the data sites. This will allow for the evaluation of the curvature and normal vector. In the MRS, a no-slip condition
\begin{equation}
\vu(\mathbf{x})=\frac{\partial \mathbf{X}}{\partial t},\label{noslip}
\end{equation}
is used to update the location of the structure. Thus, as detailed in Section 3, we need to know the velocity at some number of points along the structure. Instead of evaluating the velocity at all $N_s=50$ sample sites, we evaluate the velocity in Equation~\eqref{eq:velocity} at $N_d$ points summing over the $N_s$ sample sites. The location of the structure is then updated at the $N_d$ points $\uvXd$ using the the no slip condition given above in Equation~\eqref{noslip}. This completes a single time step. 

Representative results for a closed elastic structure immersed in an incompressible fluid with viscosity $\mu=1$ are shown in Figure \ref{fig:RSBTtest}. In the upper left panel, this is the initial configuration of the closed circle. We show results at time $t=1$, $7$, and $10$ in the other three panels. Using Equation~\eqref{eq:velocity}, we can also solve for the flow field anywhere in the 2D fluid domain. The arrows in Figure \ref{fig:RSBTtest} correspond to the flow field. The background color and colorbar corresponds to the pressure as evaluated in Equation~\eqref{eq:pressure}. Although not shown, we note that the RBF version of the MRS and the standard MRS (forces calculated using a finite difference approximation) both compare well to the SBF version of the MRS results as shown in Figure \ref{fig:RSBTtest}. 
\subsection{Stokeslet simulations on open shapes}

\subsubsection{Computation of Tangential forces}\label{sec:CortezOpen}

\begin{figure}[htb!]
\begin{center}
\includegraphics[width=2.5in]{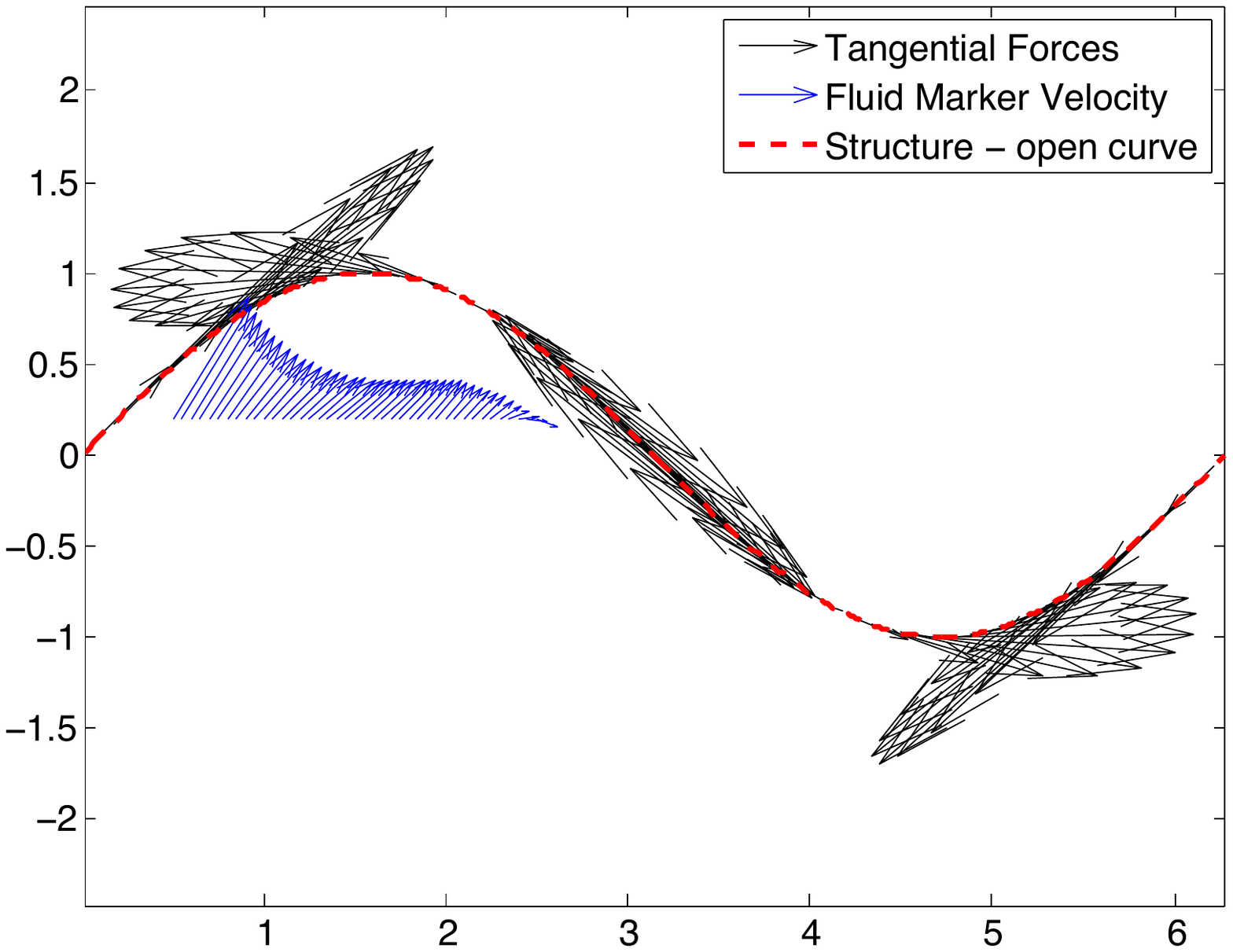}
\includegraphics[width=2.5in]{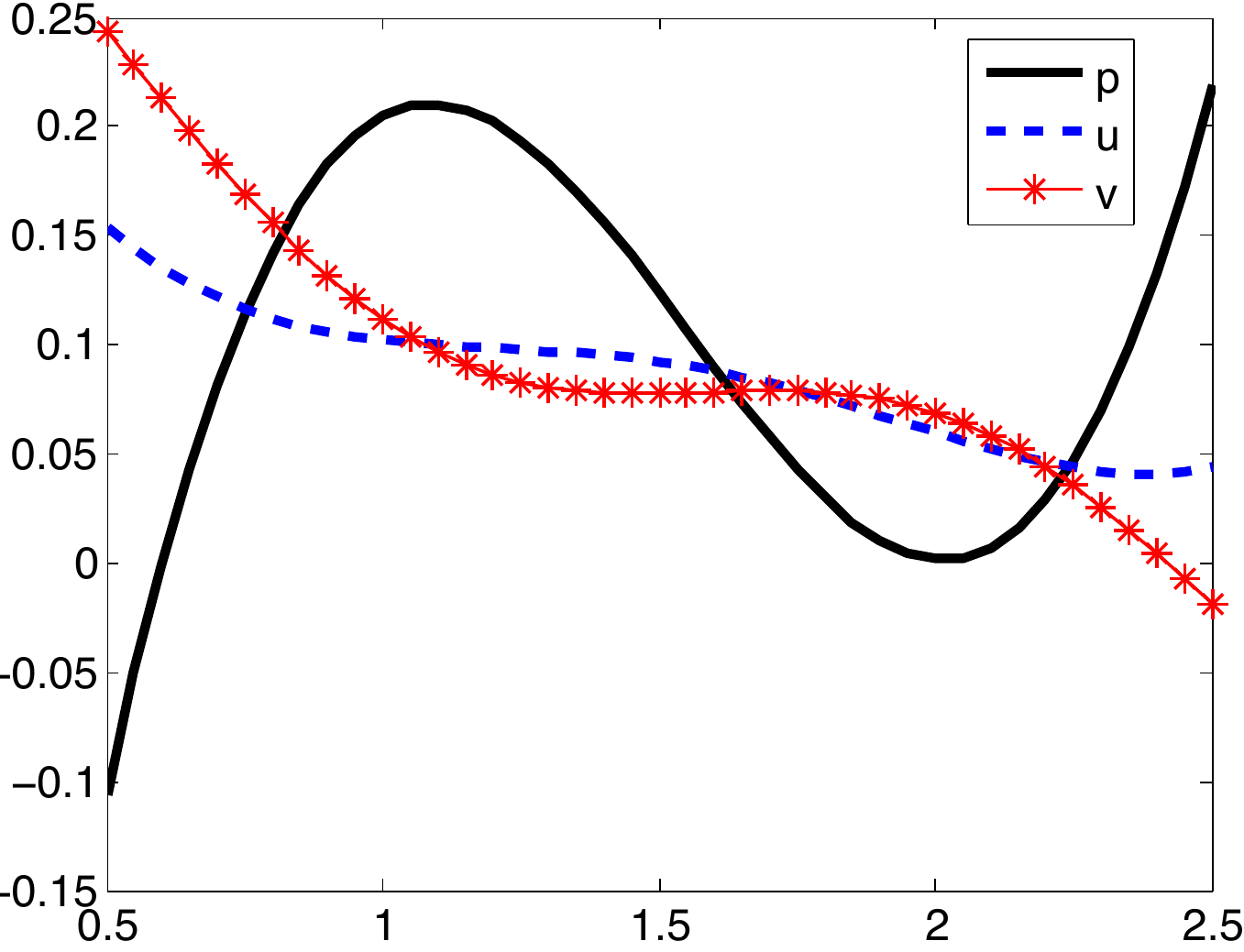}
\caption{Open curve on the left with tangential forces along the structure shown. The velocity at fluid markers along the line $y=1/5$ are also shown. On the right, exact solution for the fluid velocity at the points $y=1/5$ for the $x$ values displayed.}\label{OpenExact}
\end{center}
\end{figure}

Similar to Section \ref{sec:Cortez}, we study the same example with forces $\vF$ purely tangential to the elastic structure as defined in Equation \eqref{Fexact}.
Here, the open curve is given by $\bf{X}(\lon)=(\lon,\sin(\lon))$ for $\lon$ such that $0<\lambda<2\pi$ and we use Chebyshev and KTE nodes for the data sites. For this test case, we do not have an analytical solution for all points in the 2D fluid. However, the fundamental solution to the Stokes equation for point forces (not regularized) is given by the Stokeslet, which is defined at all points in the fluid not on the structure. We will use this as our exact solution at a set of marker points off of the structure to compare the error using different nodes. In Figure \ref{OpenExact}, the structure and tangential forces along the structure are shown. The fluid velocity is also determined along the line $y=1/5$  with viscosity set to $\mu=1$. The exact solution for the pressure and velocity along this line are shown in the right of Figure \ref{OpenExact}.
\label{sec:Results4a}

\begin{figure}
\begin{center}
\includegraphics[width=2.5in]{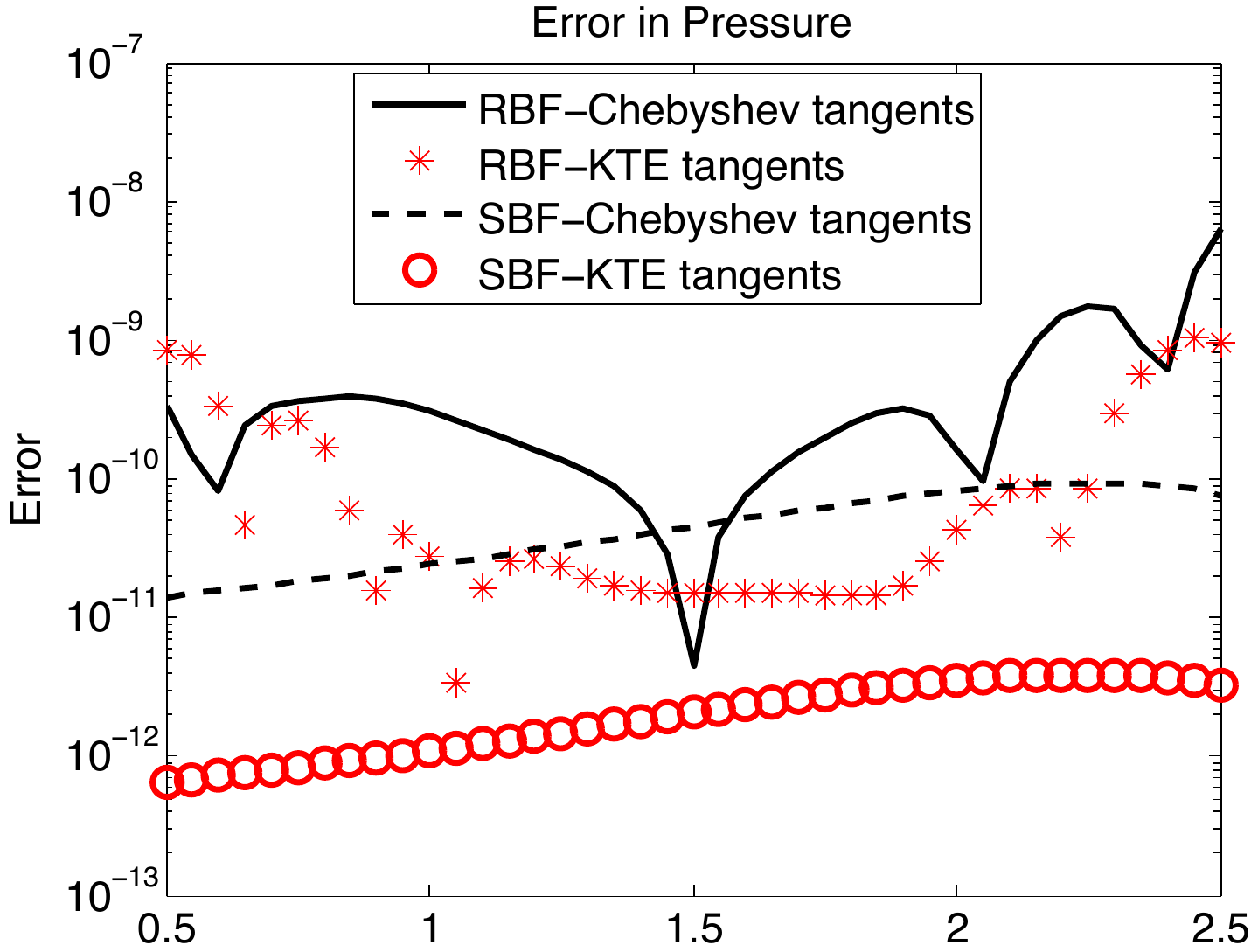}
\includegraphics[width=2.5in]{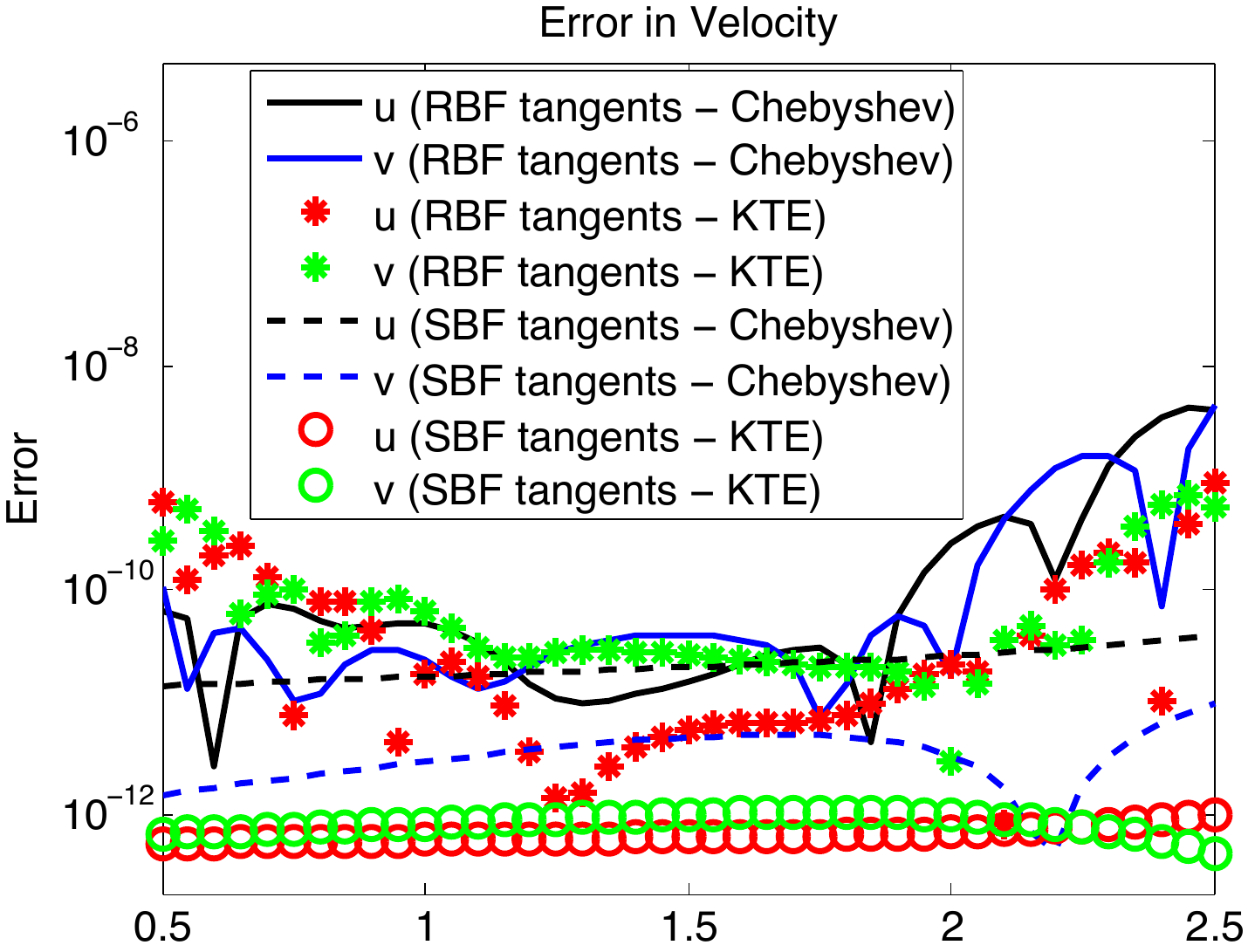}
\caption{Error of velocity at a line of marker points at $y=1/5$ for the open curve with tangential forces. The error for the pressure and velocity are shown using shape parameter $\ep=1.1$ for both Chebyshev and KTE data nodes with $\alpha=0.85$.}\label{OpenErrorExact}
\end{center}
\end{figure}

We next wanted to compare the error when evaluating the velocity and pressure when using the exact tangential forces and when using the RBF and SBF with the multiquadric kernel (for interpolation of the y-coordinate) to calculate tangential forces on the structure. In Figure \ref{OpenErrorExact}, plots of the absolute difference between the solutions at each of the x-coordinate markers are shown for pressure, $u$ component of the velocity, and $v$ component of the velocity. We use $N_d=50$ data sites (Chebyshev or KTE) to compute forces and evaluate the tangential forces at $N_s=200$ equispaced evaluation nodes. The velocity and pressure are then calculated using the fundamental solution  with forces at the evaluation nodes. The absolute error between the fundamental solution with exact tangents and RBF and SBF regularized Stokeslets solutions for the pressure and velocity are shown in Figure \ref{OpenErrorExact}. Using a shape parameter $\ep=1.1$ for both the SBF and RBF, the KTE nodes ($\alpha=0.85$) has a decreased error in the approximation of the  pressure and both components of the velocity in comparison to the Chebyshev nodes. For this $\alpha$ and $\ep$ chosen, the SBF has less error than the RBF at the maker points evaluated. Similar to the case of the closed curve in Section \ref{sec:Cortez}, the error in the pressure and velocity when using finite difference methods to calculate the force is several orders of magnitude larger when using $N_d=800$ points (results not shown). 

\subsubsection{Time-dependent simulations}
\label{sec:Results4b}
We will use a test case for an open curve corresponding to a filament propagating a planar sinusoidal wave of bending in a 2D fluid. Filaments propagating planar undulations have been studied as a simple model for a sperm flagellum \cite{Fauci95,Olson11}. The open filament will be initialized as $\vX(\lambda_k,0)=(\lambda_k,b\sin(2\pi\lambda))$ where $\lambda_k$ for $k=1,\ldots,N_s$ are KTE nodes on (0,1), $b$ is the amplitude of the sine wave, and the period of the sine wave is 1. In this test case, we want the open filament to bend or have a time dependent curvature corresponding to the sine wave $b\sin(2\pi\lambda-\omega t)$ with angular frequency $\omega$. The force on the structure will be determined as a sum of a tension force $\vF^T$ and a bending force $\vF^B$ defined as,
\begin{align}
\vF^T&=\frac{\partial}{\partial \lambda}(S_T(||\tngt||-1)\utngt),\label{Tf}\\
\vF^B&=S_B\left(\frac{\partial^4\vX}{\partial \lambda^4}-\frac{\partial \vX^I}{\partial \lambda^4}\right),\label{Bf}
\end{align}
where $\tngt$ is the tangent vector and $\utngt$ is the unit tangent. In Equation~\eqref{Tf}, this is a tension force that acts like an inextensibility constraint. The coefficients $S_T$ and $S_B$ are the tensile and bending stiffness for the elastic filament, respectively. The ideal or target shape is a propagating wave, $\bfX^I=(\lambda,b\sin(k\lambda-\omega t))$. We note that both of these forces, common to immersed boundary methods, have already been tested and implemented for a closed curve in the RBF-Immersed Boundary framework \cite{SWKFJCP2014}. Similar to the above description for the closed SBF - stokeslets method in Section \ref{closedSim} for a closed curve, we can calculate these derivatives using differentiation matrices in Equation \eqref{diffM} to determine the forces at each of the $N_s$ sample sites. We note that we will use $N_d=20$ KTE data sites with $\alpha=0.85$  and $N_s=40$ equally spaced sample sites. For this test case, we use the multiquadric radial kernel in Equation \eqref{eq:mq} and SBF given in Equation \eqref{SBFr} for interpolation of the y-coordinate.

\begin{figure}[ht]
\centering
\includegraphics[width=2.5in,height=1.5in]{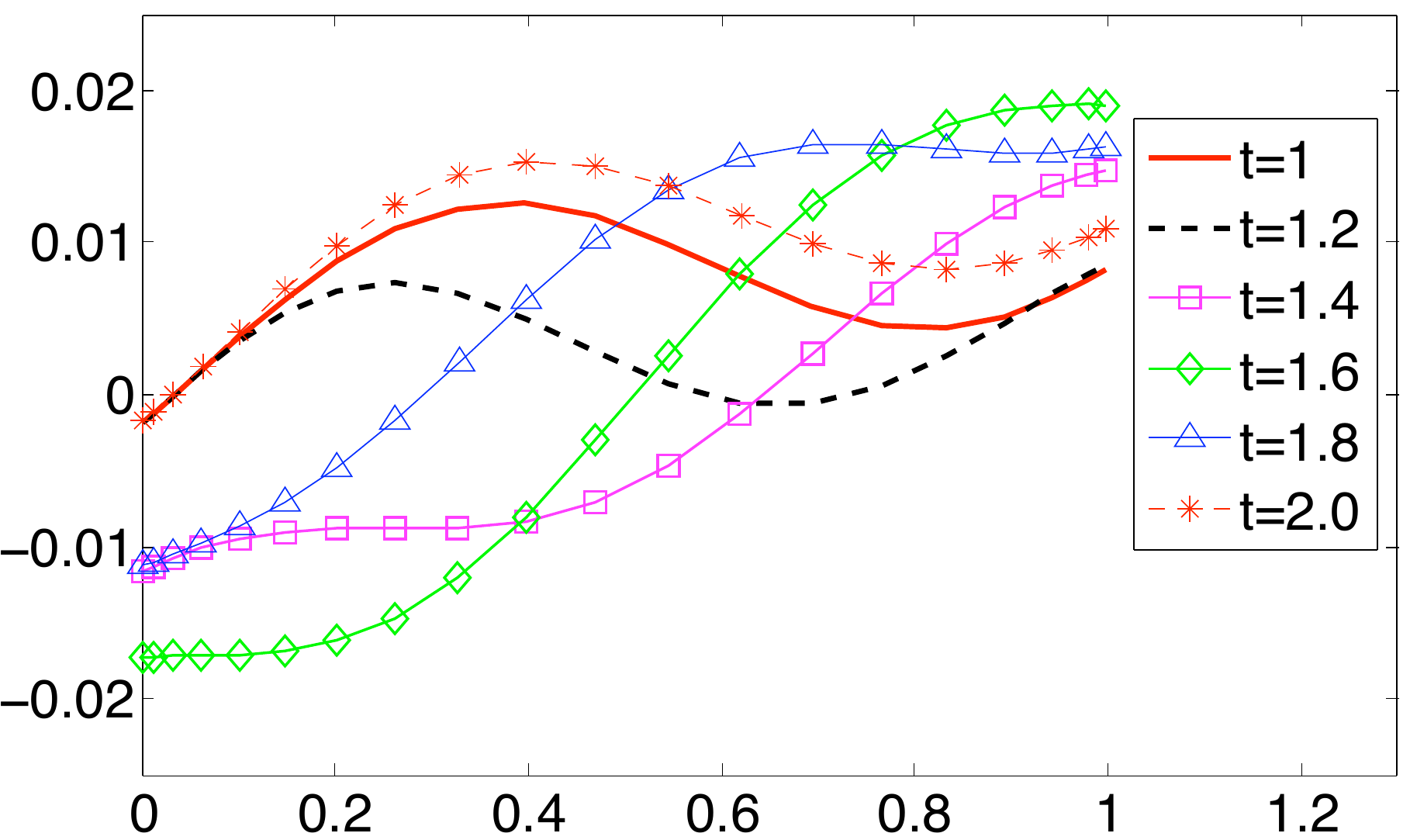}
\includegraphics[width=2.5in,height=1.5in]{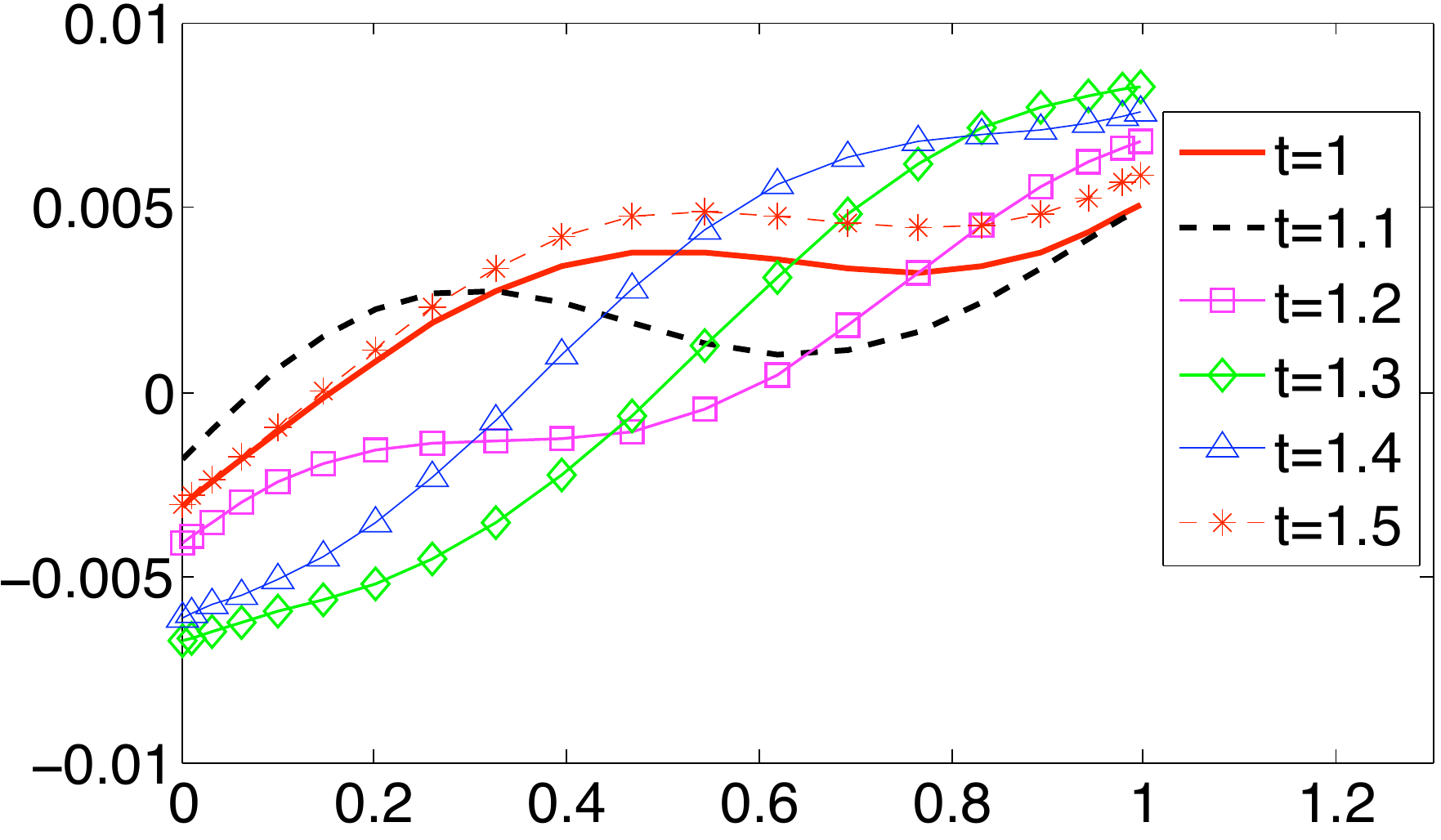}
\caption{Results for an open filament initialized with $N_d=20$, $N_s=40$, $\Delta t=5\times10^{-4}$, $\delta=2/N_s$, $S_T=0.001$, $S_B=0.1$, $\ep = 1.5$, and $\alpha=0.85$ for KTE nodes. On the left, we set $b=0.01$ and $\omega=-2\pi$. On the right, we set $b=0.005$ and $\omega=-4\pi$. }
\label{fig:OpenCaseMove}
\end{figure}

Results showing the location of the open filament at different points are shown in Figure \ref{fig:OpenCaseMove}. We solve the Stokes equations using viscosity $\mu=1$ and the no-slip condition is used to update the location of the filament. Since the bending force in Equation~\eqref{Bf} has a time dependent ideal or target location, the difference between the fourth derivative of the target location and the actual location of the structure will cause the structure to bend. On the left side of Figure \ref{fig:OpenCaseMove}, we set the amplitude of the ideal structure to be $b=0.01$ and the frequency to be $\omega=-2\pi$. As shown from the time course on the left, the filament is moving and interacting with the fluid. The forces are driving the structure to achieve an amplitude that is approximately 0.01 and it propagates a wave in one unit of time, as expected. On the right side of Figure \ref{fig:OpenCaseMove}, we have decreased the amplitude of the ideal structure to be $b=0.005$ and increased $\omega$ to -4$\pi$. As can be seen, the filament is moving with a smaller amplitude and is bending at a higher rate, corresponding to propagating a wave in 0.5 units of time. In Figure \ref{fig:OpenCaseMove}, the structure is plotted at the KTE nodes. We note that the KTE results compare well with results obtained using Chebyshev nodes with the SBF as well as using KTE or Chebyshev nodes with the RBF.

\section{Summary}

The method of Regularized Stokeslets, as a method for fluid-structure interaction at zero Reynolds number, has been successfully applied to model various biological phenomena. In this application, we presented an adaptation of the SBF geometric model developed in~\cite{Shankar2012} to the simulation of open planar curves in the method of Regularized Stokeslets, and also introduced a parametric RBF model for comparison against the SBF model. 

We conclude the following:
\begin{itemize}
\item SBF and RBF interpolants both show almost identical accuracy when used in the modeling of twice-differentiable open shapes. Given the superior accuracy of the SBF interpolant for modeling closed periodic shapes, this makes the SBF a promising choice for modeling a variety of shapes (closed and open) important to biological applications.

\item The SBF, RBF and Lagrange-Chebyshev interpolants all offer excellent accuracy and convergence in the modeling of twice-differentiable open shapes; indeed, as expected, for small numbers of points, these interpolants are superior to the Fourier interpolants that are sometimes used in the MRS (\emph{cf.},~\cite{Shankar2012}).

\item The SBF and RBF interpolants are slightly superior to the Lagrange-Chebyshev interpolant in terms of accuracy for a given number of data sites. Since SBFs and RBFs likely generalize more readily to the modeling of parametric sheets and curves embedded in $\mathbb{R}^3$, this establishes SBFs and RBFs as the more promising candidates for geometric modeling.

\item The modeling of twice-differentiable open curves requires the use of clustered points. In the case of Lagrange interpolants, the obvious choice is the set of Chebyshev nodes if the clustering is required at the ends of the interval. However, for RBF and SBF interpolants, it is clear that the mapped Chebyshev or KTE nodes are superior (for a given shape parameter). Indeed, since the KTE nodes used in our work are almost uniformly-spaced, it is reasonable to use them as general-purpose node sets for geometric modeling of both closed and open curves.

\item The dual representation afforded by global interpolants (data sites and sample sites), when used correctly, can offer greater accuracy within the method of Regularized Stokeslets for a similar or slightly lower computational cost. This result is true for all interpolants used in this work, and for Fourier interpolants as well (which are currently not used in this fashion within the MRS).

\end{itemize}

When the SBF representation was used within the IB method, a key feature was the use of both the \emph{data sites} and the \emph{sample sites} (which are analogous to IB points). The number of data sites in that work was determined by the statements on geometric accuracy made in~\cite{Shankar2012}, while the number of sample sites is chosen such that the spacing (measured in the embedding space) between any two sample sites is never more than $\frac{h}{2}$, where $h$ is the spacing of the background Eulerian grid. Since the method of regularized Stokeslets is grid-free (unlike the IB method), the number of sample sites $N_s$ is purely determined by the order of the quadrature rule used to compute the velocity (and pressure). The RBF representations (and indeed, other global representations like Fourier-based ones) will likely allow the use of higher-order quadrature rules without adversely affecting the accuracy in the computation of Lagrangian forces (since the error in the Lagrangian forces converges very rapidly with increasing $N_d$). This would also allow the use of higher-order time-stepping. We plan to explore this possibility in future work. 

We also note that error in the MRS will be due to error in the calculation of Lagrangian forces, error from discretizing the structure, or a quadrature rule for Equation~\eqref{regF}, as well as error from the regularization parameter \cite{Cortez05}. In order to have accuracy of force evaluations based on finite differences, greater accuracy is achieved provided $N_s$ is reasonably large. Thus, the accuracy of the pressure and velocity computations depend on the number of IB points $N_s$. As shown in the case of a closed curve, we get reduced error for reduced number of points when using parametric interpolants to calculate forces. It is also apparent from reviewing the RBF and MRS literature that many of the popular regularized blob functions are themselves normalized RBFs. It may be possible to exploit this connection to either enable fast evaluation of the forces, velocities and pressures, especially if the RBF used for the Lagrangian representation is the same as the RBF that serves as the blob function. We hope to explore this connection more carefully in future work.

Finally, we remark that a non-trivial extension of this current work would be to apply the SBF and RBF geometric models to the modeling of centerlines of elastic rods within the method of Regularized Stokeslets; since such modeling often requires the application of boundary conditions, the SBF and RBF geometric models need to be extended to accomplish this in a rigorous fashion. We also plan to explore this in future work.

\vspace{0.2in}
{\bf Acknowledgments:}
We would like to acknowledge useful discussions concerning this work with Professor Aaron Fogelson and the Biofluids group at the University of Utah, and with Professor Grady Wright at Boise State University. We would also like to thank the reviewers for their excellent suggestions. VS acknowledges partial support under NSF DMS-1160432. SDO acknowledges funding support under NSF DMS-1122461.

\bibliographystyle{wileyj}
\bibliography{flddoc}

\begin{thebibliography}{10}
\providecommand{\url}[1]{\texttt{#1}}
\providecommand{\urlprefix}{URL }
\expandafter\ifx\csname urlstyle\endcsname\relax
  \providecommand{\doi}[1]{doi:\discretionary{}{}{}#1}\else
  \providecommand{\doi}{doi:\discretionary{}{}{}\begingroup
  \urlstyle{rm}\Url}\fi

\bibitem{Cortez2000}
Cortez R. The method of regularized stokeslets. \emph{SIAM Journal on
  Scientific Computing}  2001; \textbf{23}(4):1204--1225,
  \doi{10.1137/S106482750038146X}.

\bibitem{Cortez05}
Cortez R, Fauci L, Medovikov A. The method of regularized {S}tokeslets in three
  dimensions: Analysis, validation, and application to helical swimming.
  \emph{Phys Fluids}  2005; \textbf{17}:0315\,041--14.

\bibitem{Cisneros08}
Cisneros L, Kessler J, Ortiz R, Cortez R, Bees M. Unexpected bipolar flagellar
  arrangements and long-range flows driven by bacteria near solid boundaries.
  \emph{Phys Rev Lett}  2008; \textbf{101}:168\,102--1--4.

\bibitem{Flores05}
Flores H, Lobaton E, Mendez-Diez S, Tlupova S, Cortez R. A study of bacterial
  flagellar bundling. \emph{Bull Math Bio}  2005; \textbf{65}:137--168.

\bibitem{Jung07}
Jung S, Mareck K, Fauci L, Shelley M. Rotational dynamics of a superhelix towed
  in a stokes fluid. \emph{Phys Fluids}  2007; \textbf{19}:103\,105--1--6.

\bibitem{Nguyen11b}
Nguyen H, Ortiz R, Cortez R, Fauci L. The action of waving cylindrical rings in
  a viscous fluid. \emph{J Fluid Mech}  2011; \textbf{671}:574--586.

\bibitem{Gillies09}
Gillies E, Cannon R, Green R, Pacey A. Hydrodynamic propulsion of human sperm.
  \emph{J Fluid Mech}  2009; \textbf{625}:445--474.

\bibitem{Olson11}
Olson S, Suarez S, Fauci L. Coupling biochemistry and hydrodynamics captures
  hyperactivated sperm motility in a simple flagellar model. \emph{J Theor Bio}
   2011; \textbf{283}:203--216.

\bibitem{Olson13a}
Olson S, Lim S, Cortez R. Modeling the dynamics of an elastic rod with
  intrinsic curvature and twist using a regularized stokes formulation. \emph{J
  Comp Phys}  2013; \textbf{283}:69--187.

\bibitem{Olson13b}
Olson S. Fluid dynamic model of invertebrate sperm motility with varying
  calcium inputs. \emph{J Biomechanics}  2013; \textbf{46}:329--337.

\bibitem{Olson14a}
Olson S. Motion of filaments with planar and helical bending waves in a viscous
  fluid. \emph{Biological Fluid Dynamics: Modeling, Computation, and
  Applications, AMS Contemp Math}, Layton A, Olson S (eds.), AMS: Providence,
  RI, 2014; ???

\bibitem{Simons14}
Simons J, Olson S, Cortez R, Fauci L. The dynamics of sperm detachment from
  epithelium in a coupled fluid-biochemical model of hyperactivated motility.
  \emph{J Theor Biol}  2014; \textbf{354}:81--94.

\bibitem{Ainley08}
Ainley J, Durkin S, Embid R, Boindala P, Cortez R. The method of images for
  regularized {S}tokeslets. \emph{J Comp Phys}  2008; \textbf{227}:4600--4616.

\bibitem{Leiderman13}
Leiderman K, Bouzarth EL, Cortez R, Layton AT. A regularization method for the
  numerical solution of periodic {S}tokes flow. \emph{J Comput Phys}  2013;
  \textbf{236}:187--202.

\bibitem{Leiderman14}
Leiderman K, Bouzarth EL, Nguyen H. A regularization method for the numerical
  solution of doubly-periodic stokes flow. \emph{Biological Fluid Dynamics:
  Modeling, Computation, and Applications, AMS Contemp Math}, Layton A, Olson S
  (eds.), AMS: Providence, RI, 2013; ???

\bibitem{Cogan05}
Cogan NG, Cortez R, Fauci LJ. Modeling physiological resistance in bacterial
  biofilms. \emph{Bull. Math. Biol.}  2005; \textbf{67}:831--853.

\bibitem{Bouzarth11}
Bouzarth EL, Layton AT, Young YN. Modeling a semi-flexible filament in cellular
  {S}tokes flow using regularized {S}tokeslets. \emph{Int J Numer Methods
  Biomed Eng}  2012; \textbf{12}:2021--2034.

\bibitem{Lee14}
Lee W, Kim Y, Olson S, S L. Nonlinear dynamics of a rotating elastic in a
  viscous fluid. \emph{Phys Rev E}  2014; \textbf{90}:033\,012.

\bibitem{OMalley12}
O'Malley S, Bees M. The orientation of swimming biflagellates in shear flows.
  \emph{Bull Math Bio}  2012; \textbf{74}:232--255.

\bibitem{Shankar2012}
Shankar V, Wright GB, Fogelson AL, Kirby RM. A study of different modeling
  choices for simulating platelets within the immersed boundary method.
  \emph{Applied Numerical Mathematics}  2013; \textbf{63}(0):58 -- 77,
  \doi{10.1016/j.apnum.2012.09.006}.

\bibitem{SWKFJCP2014}
Shankar V, Wright GB, Kirby RM, Fogelson AL. Augmenting the immersed boundary
  method with radial basis functions (rbfs) for the modeling of platelets in
  hemodynamic flows. Submitted.

\bibitem{DriscollFornberg2002}
Driscoll T, Fornberg B. Interpolation in the limit of increasingly flat radial
  basis functions. \emph{Computers \& Mathematics with Applications}  2002;
  \textbf{43}(3):413--422.

\bibitem{FasshauerSchumaker:1998}
Fasshauer GE, Schumaker LL. Scattered data fitting on the sphere.
  \emph{Mathematical Methods for Curves and Surface, Vol.2 of the Proceedings
  of the 4th Int. Conf. on Mathematical Methods for Curves and Surfaces,
  Lillehammer, Norway}, Daehlen M, Lyche T, Schumaker LL (eds.), Vanderbilt
  University Press, Nashville Tennessee, 1998.

\bibitem{Wendland:2004}
Wendland H. \emph{Scattered data approximation}, \emph{Cambridge Monographs on
  Applied and Computational Mathematics}, vol.~17. Cambridge University Press:
  Cambridge, 2005.

\bibitem{JetterStocklerWard:1999}
Jetter K, St\"{o}ckler J, Ward JD. Error estimates for scattered data
  interpolation on spheres. \emph{Math. Comput.}  April 1999;
  \textbf{68}(226):733--747.

\bibitem{NarcSunWard:2007}
Narcowich FJ, Sun X, Ward JD, Wendland H. Direct and inverse {S}obolev error
  estimates for scattered data interpolation via spherical basis functions.
  \emph{Found. Comput. Math.}  2007; \textbf{7}(3):369--390.

\bibitem{FornbergPiret:2007}
Fornberg B, Piret C. A stable algorithm for flat radial basis functions on a
  sphere. \emph{SIAM J. Sci. Comp.}  2007; \textbf{30}:60--80.

\bibitem{FuselierWright2014}
Fuselier E, Wright G. Order-preserving derivative approximation with periodic
  radial basis functions. \emph{Advances in Computational Mathematics}  2014;
  :1--31\doi{10.1007/s10444-014-9348-1}.
  \urlprefix\url{http://dx.doi.org/10.1007/s10444-014-9348-1}.

\bibitem{SWFKIJNMF2014}
Shankar V, Wright GB, Fogelson AL, Kirby RM. A radial basis function (rbf)
  finite difference method for the simulation of reaction–diffusion equations
  on stationary platelets within the augmented forcing method.
  \emph{International Journal for Numerical Methods in Fluids}  2014;
  \textbf{75}(1):1--22, \doi{10.1002/fld.3880}.
  \urlprefix\url{http://dx.doi.org/10.1002/fld.3880}.

\bibitem{FuselierWright2012}
Fuselier EJ, Wright GB. A high-order kernel method for diffusion and
  reaction-diffusion equations on surfaces. \emph{Journal of Scientific
  Computing}  2013; :1--31\doi{10.1007/s10915-013-9688-x}.
  \urlprefix\url{http://dx.doi.org/10.1007/s10915-013-9688-x}.

\bibitem{SWFKJSC2014}
Shankar V, Wright GB, Kirby RM, Fogelson AL. A radial basis function
  (rbf)-finite difference (fd) method for diffusion and reaction-diffusion
  equations on surfaces. Submitted.

\bibitem{Fasshauer:2007}
Fasshauer GE. \emph{Meshfree Approximation Methods with {MATLAB}}.
  Interdisciplinary Mathematical Sciences - Vol. 6, World Scientific
  Publishers: Singapore, 2007.

\bibitem{Peskin:2002}
Peskin CS. The immersed boundary method. \emph{Acta Numerica}  2002;
  \textbf{11}:479--517.

\bibitem{PlatteAnalytic}
Platte RB. How fast do radial basis function interpolants of analytic functions
  converge? \emph{IMA Journal of Numerical Analysis}  2011;
  \textbf{31}(4):1578--1597.

\bibitem{OlsonLimCortezJCP2013}
Olson SD, Lim S, Cortez R. Modeling the dynamics of an elastic rod with
  intrinsic curvature and twist using a regularized stokes formulation.
  \emph{Journal of Computational Physics}  2013; \textbf{238}(0):169 -- 187,
  \doi{http://dx.doi.org/10.1016/j.jcp.2012.12.026}.

\bibitem{Cortez00}
Cortez R, Minion M. The blob projection method for immersed boundary problems.
  \emph{J Comput Phys}  2000; \textbf{161}:428--453.

\bibitem{Cortez12}
Cortez R, Nicholas M. Slender body theory for stokes flows with regularized
  forces. \emph{CAMCoS}  2012; \textbf{7}:33--62.

\bibitem{Fauci95}
Fauci L, Mc{D}onald A. Sperm motility in the presence of boundaries. \emph{Bull
  Math Biol}  1995; \textbf{57}:679--699.

\end{thebibliography}
 
\end{document}